\documentclass[11pt]{article}

\usepackage{latexsym}
\usepackage{graphicx}
\usepackage{framed}
\usepackage{threeparttable}
\usepackage{amssymb}
\usepackage{amsmath}
\usepackage{multirow}
\usepackage{hyperref}
\usepackage{tikz}

\begin{document}

\newcommand{\noi}{\noindent}
\newcommand{\nn}{\nonumber}
\newcommand{\bd}{\begin{displaymath}}
\newcommand{\ed}{\end{displaymath}}
\newcommand{\bp}{\underline{\bf Proof}:\ }
\newcommand{\ep}{{\hfill $\Box$}\\ }
\newtheorem{1}{LEMMA}[section]
\newtheorem{2}{THEOREM}[section]
\newtheorem{3}{COROLLARY}[section]
\newtheorem{4}{PROPOSITION}[section]
\newtheorem{5}{REMARK}[section]
\newtheorem{20}{OBSERVATION}[section]
\newtheorem{10}{DEFINITION}[section]
\newtheorem{30}{RESULTS}[section]
\newtheorem{40}{CLAIM}[section]
\newtheorem{50}{ASSUMPTION}[section]
\newtheorem{60}{EXAMPLE}[section]
\newtheorem{70}{ALGORITHM}[section]
\newtheorem{80}{PROBLEM}
\newcommand{\be}{\begin{equation}}
\newcommand{\ee}{\end{equation}}
\newcommand{\ba}{\begin{array}}
\newcommand{\ea}{\end{array}}
\newcommand{\bea}{\begin{eqnarray}}
\newcommand{\eea}{\end{eqnarray}}
\newcommand{\bqn}{\begin{eqnarray*}}
\newcommand{\eqn}{\end{eqnarray*}}


\newcommand{\e} { \ = \ }
\newcommand{\leqs}{ \ \leq \ }
\newcommand{\geqs}{ \ \geq \ }
\def\theequation{\thesection.\arabic{equation}}
\def\bReff#1{{\bRm
(\bRef{#1})}}
\newcommand{\eps}{\varepsilon}
\newcommand{\sgn}{\operatorname{sgn}}
\newcommand{\sign}{\operatorname{sign}}
\newcommand{\Vol}{\operatorname{Vol}}
\newcommand{\Var}{\operatorname{Var}}
\newcommand{\Cov}{\operatorname{Cov}}
\newcommand{\vol}{\operatorname{vol}}
\newcommand{\var}{\operatorname{var}}
\newcommand{\cov}{\operatorname{cov}}
\renewcommand{\Re}{\operatorname{Re}}
\renewcommand{\Im}{\operatorname{Im}}
\newcommand{\bE}{{\mathbb E}}
\newcommand{\bR}{\mathbb{R}}
\newcommand{\bN}{{\mathbb N}}
\newcommand{\bC}{\mathbb{C}}
\newcommand{\bF}{\mathbb{F}}
\newcommand{\bQ}{{\mathbb Q}}
\newcommand{\bZ}{{\mathbb Z}}
\newcommand{\cA}{{\mathcal A}}
\newcommand{\cB}{{\mathcal B}}
\newcommand{\cC}{{\mathcal C}}
\newcommand{\cD}{{\mathcal D}}
\newcommand{\cE}{{\mathcal E}}
\newcommand{\cF}{{\mathcal F}}
\newcommand{\cG}{{\mathcal G}}
\newcommand{\cH}{{\mathcal H}}
\newcommand{\cI}{{\mathcal I}}
\newcommand{\cJ}{{\mathcal J}}
\newcommand{\cK}{{\mathcal K}}
\newcommand{\cL}{{\mathcal L}}
\newcommand{\cM}{{\mathcal M}}
\newcommand{\cN}{{\mathcal N}}
\newcommand{\cO}{{\mathcal O}}
\newcommand{\cP}{{\mathcal P}}
\newcommand{\cQ}{{\mathcal Q}}
\newcommand{\cR}{{\mathcal R}}
\newcommand{\cS}{{\mathcal S}}
\newcommand{\cT}{{\mathcal T}}
\newcommand{\cU}{{\mathcal U}}
\newcommand{\cV}{{\mathcal V}}
\newcommand{\cW}{{\mathcal W}}
\newcommand{\cX}{{\mathcal X}}
\newcommand{\cY}{{\mathcal Y}}
\newcommand{\cZ}{{\mathcal Z}}
\newcommand{\bx}{{\mathbf x}}
\newcommand{\by}{{\mathbf y}}
\newcommand{\bz}{{\mathbf z}}
\newcommand{\bba}{{\mathbf a}}
\newcommand{\bbb}{{\mathbf b}}
\newcommand{\bbc}{{\mathbf c}}


\title{Computing tensor eigenvalues via homotopy methods}
\author{ Liping Chen\thanks{Department of Mathematics, Michigan State University, East Lansing, MI 48824, USA. Email: \texttt{chenlipi@msu.edu}}  \and
Lixing Han\thanks{Department of Mathematics, University of Michigan-Flint, Flint, MI 48502, USA. Email: \texttt{lxhan@umflint.edu}} \and Liangmin Zhou\thanks{Department of Mathematics, Michigan State University, East Lansing, MI 48824, USA. Email: \texttt{zhoulian@msu.edu}}
 }
\date{}
\maketitle

\begin{abstract} 

We introduce the concept of mode-$k$ generalized eigenvalues and eigenvectors of a tensor and prove some properties of such eigenpairs. In particular, we derive an upper bound for the number of equivalence classes of generalized tensor eigenpairs using mixed volume. Based on this bound and the structures of tensor eigenvalue problems, we propose two homotopy continuation type algorithms to solve tensor eigenproblems. With proper implementation, these methods can find all equivalence classes of isolated generalized eigenpairs and some generalized eigenpairs contained in the positive dimensional components (if there are any). We also introduce an algorithm that combines a heuristic approach and a Newton homotopy method to extract real generalized eigenpairs from the found complex generalized eigenpairs. A MATLAB software package \texttt{TenEig} has been developed to implement these methods. Numerical results are presented to illustrate the effectiveness and efficiency of \texttt{TenEig} for computing complex or real generalized eigenpairs.

\end{abstract}



\section{Introduction}
\label{Intro}
\setcounter{equation}{0}

Eigenvalues of tensors were first introduced by   Lim \cite{Lim05} and Qi \cite{Qi05} in 2005. Since then, tensor eigenvalues have found applications in automatic control, statistical data analysis,  diffusion tensor imaging,  image authenticity verification, spectral hypergraph theory, and quantum entanglement,   etc., see for example, \cite{CQZ13, CD12, HQZ12, Qi05, QSW07, QWW08, QYW10, QYX} and the references therein.  The tensor eigenvalue problem has become an important subject of numerical multilinear algebra.  

Various definitions of eigenvalues for tensors have been proposed in the literature,  including E-eigenvalues  and eigenvalues  in the complex field, and Z-eigenvalues, H-eigenvalues, and D-eigenvalues in the real field \cite{Lim05, Qi05, QWW08}.  In \cite{CPZ09}, Chang,  Pearson, and Zhang introduced a notion of generalized eigenvalues for tensors that unifies several types of eigenvalues. Recently this definition has been further generalized by Cui, Dai, and Nie \cite{CDN14}.

Unlike the matrix eigenvalue problem, computing eigenvalues of the third or higher order tensors  is a difficult problem \cite{HL13}. Nonetheless, several algorithms which aim at computing one or some eigenvalues of a tensor have been developed recently. These algorithms are designed for tensors of certain type, such as entry-wise nonnegative or symmetric tensors.   

For nonnegative tensors, Ng, Qi, and Zhou \cite{NQZ09} proposed a power-type method for computing the largest H-eigenvalue of a nonnegative tensor. Modified versions of the Ni-Qi-Zhou method have been proposed in \cite{LZI10, ZQX12, ZQW13}.  

For real symmetric tensors, Hu, Huang, and Qi \cite{HHQ13} proposed a sequential semidefinite programming method for computing extreme Z-eigenvalues. 
 Hao, Cui, and Dai \cite{HCD15} proposed a sequential subspace projection method for a similar purpose. 
Kolda and Mayo \cite{KM11} proposed a shifted power method (SSHOPM) for computing a Z-eigenvalue.  They have improved SSHOPM in \cite {KM14} by updating the shift parameter adaptively. The resulting method can be used to compute a real generalized eigenvalue.   Han \cite{Han13} proposed an unconstrained optimization method for computing a real generalized eigenvalue for even order real symmetric tensors. The methods in \cite{Han13, KM11, KM14} can find more eigenvalues of a symmetric tensor if they are run multiple times using different starting points. Recently,  Cui, Dai, and Nie \cite{CDN14} proposed a novel method for computing all real generalized eigenvalues. 

In this paper, we are concerned with computing all eigenpairs of a general real or complex tensor. As can be seen from the next section, finding eigenpairs of a tensor amounts to solving a system of polynomials. Naturally one would consider to use algebraic geometry methods  such as the Gr\"{o}bner basis method and the resultant method \cite{CLO07} for this purpose. These methods can obtain symbolic solutions of a polynomial system, which are accurate. However, they are expensive in terms of computational cost and space. Moreover, they are difficult to parallelize.   A class of numerical methods, the homotopy continuation methods, can overcome these shortcomings of the Gr\"{o}bner basis and the resultant methods.   During the past few decades, significant advances have been made on homotopy continuation methods for polynomial systems, see for example, \cite{BHSW,LLT08, Li03, Morgan09, SW05}. Recently, the homotopy techniques have been used to study tensor decomposition and perfect identification problems (\cite{HOOS15}).    
  
In this paper we investigate  computing complex eigenpairs of general tensors using homotopy continuation methods.  One attractive feature of the homotopy continuation methods is that they can find all isolated solutions of polynomial systems and  some solutions in the positive dimensional solution components. We propose two homotopy-type algorithms  for computing complex eigenpairs of a tensor. These algorithms allow us to find all equivalence classes of isolated eigenpairs of a general tensor and  some eigenpairs in positive dimensional eigenspaces (if there are any). We also present an algortihm combining a heuristic approach and a Newton homotopy method
to compute real eigenpairs based on the found complex eigenpairs. Numerical examples show that our methods are effective and efficient.

This paper is organized as follows. In Section 2, we define mode-$k$ generalized eigenvalues and eigenvectors which extend the matrix right eigenpairs and left eigenpairs to higher order tensors. Some properties of such eigenpairs are proved. An upper bound for the number of equivalence classes of generalized tensor eigenpairs using mixed volume is derived. In Section 3, we consider computing mode-$k$ generalized complex eigenpairs and present two algorithms.  In Section 4, we introduce a method  to compute real mode-$k$ generalized eigenpairs. Finally in Section 5, some numerical results are provided.

\section{Tensor eigenvalues and eigenvectors}
\label{sec2}
\setcounter{equation}{0}

 Let $\bF = \bC \ {\rm or} \ \bR $  be the complex field or the real field.  Let $m \geq 2$, $m' \geq 2$, and $n$ be positive integers. Denote the set of all $m$th-order, $n$-dimensional  tensors on the field $\bF$ by $\bF^{[m,n]}$. 
A tensor in    $\bF^{[m,n]}$ is indexed as
$$ 
\cA = (A_{i_1 i_2 \cdots i_m}), 
$$  
where $ A_{i_1 i_2 \cdots i_m} \in \bF$, for $ 1 \leq i_1, i_2, \cdots, i_m \leq n$. 

For $x \in \bC^n$, the tensor $\cA$ defines a multilinear form
\be
\label{form}
 \cA x^m = \sum_{i_1, \cdots, i_m =1}^{n} A_{i_1 i_2 \cdots i_m}x_{i_1} x_{i_2} \cdots x_{i_m}. 
\ee
For $1 \leq k \leq m$, $\cA^{(k)} x^{m-1}$ is an $n$-vector whose $j$th entry is defined as
\be
\label{modekvec}
(\cA^{(k)} x^{m-1})_j = \sum_{i_1, \cdots, i_{k-1},  i_{k+1}, \cdots, i_m =1}^{n} A_{i_1  \cdots i_{k-1}  j i_{k+1} \cdots  i_m}x_{i_1} \cdots x_{i_{k-1}}  x_{i_{k+1}} \cdots x_{i_m}.
\ee
When $k=1$, the vector $\cA^{(1)} x^{m-1} $ is denoted by $\cA x^{m-1} $. 

A  real tensor $\cA \in \bR^{[m,n]}$ is {\it positive definite}  if the multilinear form 
$\cA x^m $
is positive for all $x \in \bR^n \backslash \{ 0\}$. A tensor $\cA \in \bF^{[m,n]}$ is {\it symmetric} if its entries $A_{i_1 i_2 \cdots i_m} $ are invariant under any permutations of their indices $i_1, i_2, \cdots, i_m$. 

We now introduce  the following mode-$k$ generalized eigenvalue definition for a general tensor $\cA$.

\begin{10}  
\label{def1}
  Let $\cA \in \bF^{[m,n]}$ and $\cB \in  \bF^{[m',n]}$.    Assume that  $ \cB x^{m'}$ is not identically zero as a function of $x$. For $ 1 \leq k \leq m$, 
  if there exist a scalar $\lambda \in \bC$ and a vector $ x  \in \bC^n \backslash \{0\} $ such that
\begin{itemize}
\item when $m=m'$,
\be
\label{tevalue}
\cA^{(k)} x^{m-1} = \lambda \cB x^{m-1}, 
\ee
\item when $m \ne m'$,
\be
\label{evalue}
\cA^{(k)} x^{m-1} = \lambda \cB x^{m'-1}, \ \ \ \cB x^{m'} = 1,
\ee
\end{itemize}
 then $\lambda$ is called a mode-$k$ $\cB$-eigenvalue of $\cA$ and $x$ a mode-$k$ $\cB$-eigenvector associated with $\lambda$.  $(\lambda, x)$ is called a mode-$k$ $\cB$-eigenpair of $\cA$. 

If  
$\lambda \in \bR, x \in \bR^n$, then $\lambda$ is called a mode-$k$ $\cB_R$-eigenvalue of $\cA$ and $x$ a mode-$k$ $\cB_R$-eigenvector associated with $\lambda$, and  $(\lambda, x)$  a mode-$k$ $\cB_R$-eigenpair of $\cA$.

Denote the set of all  mode-$k$ $\cB$ eigenvalues of $\cA$ by $\sigma_{\cB}({\cA^{(k)}})$. 
\end{10}

\begin{5}
\label{remark4} 
{\rm 
Let $(\lambda,x)$ be a mode-$k$ $\cB$-eigenpair of $\cA $. By (\ref{tevalue}) or (\ref{evalue}), $(\lambda,x)$ is a solution to  $\cA^{(k)} x^{m-1}=\lambda \cB x^{m'-1}$. So is $(\lambda',x')$ with $\lambda'=t^{m-m'}\lambda$ and $x'=tx$ for $t\in \bC \backslash \{0\}$. From this point of view, the solution space of $\cA^{(k)} x^{m-1}=\lambda \cB x^{m'-1}$ consists of different equivalence classes.  We denote such an equivalence class by
 $$[(\lambda,x)]:=\{(\lambda',x') \, | \, \lambda'=t^{m-m'}\lambda, x'=tx, t \in \bC \backslash \{0\}\}.$$ 

When $m \ne m'$, taking arbitrary $(\lambda',x')\in [(\lambda,x)]$ and substituting $x'=tx$ into $\cB x^{m'}=1$ in (\ref{evalue}) yields $t^{m'}=1$, which gives $m'$ different values for $t$. This implies that the normalization $\cB x^{m'}=1$ in (\ref{evalue}) restricts us to choose $m'$ representative solutions from each equivalence class. 

In our later discussions, we often choose only one representative from each equivalence class, and we often count the number of equivalence classes of mode-$k$ $\cB$-eigenpairs.
}
\end{5}

\begin{5}
\label{remark5} 
{\rm 
If only one representative is desirable from each equivalence class of eigenpairs, we can solve  $\cA^{(k)} x^{m-1}=\lambda \cB x^{m'-1}$ augmented with an additional linear equation
\be
\label{hyperplane} 
a_1 x_1 + a_2x_2+ \cdots + a_n x_n+b =0,
\ee
where $a_1, \dots, a_n,b$ are random complex numbers. Then normalize the resulting solution to satisfy $\cB x^{m'}=1$ in the case $m \ne m'$.
}
\end{5}

In the matrix case when $m=m'=2$ and $\cB=I_n$ (the $n \times n$ identity matrix), the mode-1 eigenvectors are right eigenvectors and the mode-2 eigenvectors are left eigenvectors of $\cA$, and the mode-1 and mode-2 eigenvalues are the eigenvalues of matrix $\cA$, i.e., $\sigma_{\cB}({\cA^{(1)}}) = \sigma_{\cB}({\cA^{(2)}})$. However, when $m \geq 3$, $\sigma_{\cB}({\cA^{(k)}}) = \sigma_{\cB}({\cA^{(l)}})$ is generally not true when $k \ne l$, unless $\cA$ has a certain type of symmetry. The following example illustrates this situation. 

\begin{60}
\label{example1}
{\rm
Consider the tensor $\cA \in \bR^{[3,2]}$ whose entries are 
\[
\begin{array}{c} 
A_{111}= 1, A_{121}=  2, A_{211}= 3, A_{221}= 4, \\
A_{112}= 5, A_{122}=   6, A_{212}= 7, A_{222}=  0.
\end{array}
\]
Choose $m'=2$ and $\cB=I_2$ (the $2 \times 2$ identity matrix). Note that in this case, if $(\lambda, x)$ is an $\cB$-eigenpair of $\cA$, so is $(-\lambda, -x)$. We follow \cite{CS13}, regarding $(\lambda, x)$ and $(-\lambda, -x)$ as the  same eigenpair. Then
$$
\sigma_{\cB} (\cA^{(1)}) = \{
   0.4105, 
   4.3820,
   9.8995 \},
$$
$$
\sigma_{\cB} (\cA^{(2)}) = \{
   0.2851,                  
   4.3536,
   9.5652
\},
$$
$$
\sigma_{\cB} (\cA^{(3)}) = \{
   0.2936,
   4.3007,
   9.4025
\}.
$$
Clearly, $\sigma_{\cB}({\cA^{(k)}}) \ne \sigma_{\cB}({\cA^{(l)}})$  when $k \ne l$. 
}
\end{60}

\begin{4}
\label{proposition1}
 Suppose that $(\lambda, x)$ is a mode-$k$ $\cB$-eigenpair, $(\mu, x)$ is a mode-$l$ $\cB$-eigenpair of $\cA$, and when $m=m'$, $\cB x^m \ne 0$. Then $\lambda = \mu$.  
\end{4}
\bp
Note that 
$$
\lambda \cB x^{m'} = \cA x^m = \mu \cB x^{m'}.
$$
This immediately implies that $\lambda = \mu$, since $\cB x^{m'}=1$ when $m \ne m'$ and $\cB x^m \ne 0$ when $m=m'$. 
\ep

Let $\cA \in \bF^{[m,n]}$. For $1 \leq k < l \leq m$, tensor $\cG \in \bF^{[m,n]}$ 
is said to be the $\langle k, l \rangle$ {\it transpose} of $\cA$ if 
$$
\cG_{i_1 \cdots i_{k-1}  i_l  i_{k+1} \cdots i_{l-1}  i_k  i_{l+1}  \cdots i_m} = \cA_{i_1  \cdots i_{k-1}  i_k i_{k+1} \cdots i_{l-1} i_l  i_{l+1} \cdots i_m},
$$
for all $ 1\leq i_1, \cdots, i_m \leq m$. Denote the $\langle k,l \rangle$ {\it transpose} of $\cA$ by $\cA^{\langle k,l \rangle}$. We say that tensor $\cA$ is $\langle k,l \rangle$ {\it partially symmetric} if 
$$
\cA^{\langle k,l \rangle} = \cA.
$$

\begin{4}
\label{theorem1}
 Let $\cA \in \bF^{[m,n]}$ and $\cB \in  \bF^{[m',n]}$.    Assume that  $\cB x^{m'}$ is not identically zero as a function of $x$.  Let $k,l$ be integers such that $1 \leq k < l \leq m$. Then
\begin{itemize}
\item  $(\lambda, x)$ is a mode-$k$ $\cB$-eigenpair of $\cA$ if and only if it is a  mode-$l$ $\cB$-eigenpair of $\cA^{\langle k,l \rangle}$. 
\item The sets of mode-$k$ $\cB$-eigenpairs and mode-$l$ $\cB$-eigenpairs are the same if $\cA$ is $\langle k,l \rangle$  partially symmetric. 
\end{itemize}
\end{4}

The  eigenvalues/eigenvectors defined in  \cite{CPZ09, CDN14, Qi05, QWW08} are mode-1 eigenvalues/eigenvectors. The tensors considered in these papers are primarily real symmetric tensors. For symmetric tensors, the sets of mode-$k$ $\cB$-eigenpairs and  mode-$1$ $\cB$-eigenpairs are the same for any $k$. Therefore, mode-1 eigenvalues serve the purpose of those papers. On the other hand, nonsymmetric tensors arise from applications and theoretical studies, see, for example, \cite{CS13, CPZ08, FGH13, NQZ09, YY10, YY11}.   In \cite{Lim05}, Lim  defined mode-$k$ eigenvalues/eigenvectors for nonsymmetric real tensors $\cA$ when $\cB $ is the $m'$th order identity tensor for some $m' \geq 2$. Definition \ref{def1} considers more general $\cA$ and $\cB$.

As in \cite{CPZ09, CDN14}, Definition \ref{def1} adapts a unified approach to define tensor eigenvalues. It covers various types of tensor eigenvalues introduced in the literature, including 
  
\begin{itemize}

\item If $\cA \in \bR^{[m,n]}$, $m'=2$, and $\cB$ is the identity matrix $I_n \in \bR^{n\times n} $, the mode-1 $\cB$-eigenpairs are the E-eigenpairs and the mode-1 $\cB_R$-eigenpairs are the Z-eigenpairs defined in  \cite{Qi05}, which satisfy
\be
\label{Zeig}
\cA x^{m-1} = \lambda x, \ \ \ x^Tx=1.
\ee

\item If $\cA \in \bR^{[m,n]}$, $m'=2$ and $\cB = D$, where $D \in \bR^{n\times n} $ is  a symmetric positive definite  matrix, 
the $\cB_R$-eigenpairs are the D-eigenpairs defined in \cite{QWW08}, which satisfy 
\be
\label{Deig}
\cA x^{m-1} = \lambda Dx, \ \ \ x^TDx=1.
\ee  

\item  If $\cA \in \bR^{[m,n]}$, $m=m'$ and $\cB =\cI$  is the identity tensor,  mode-1 $\cB$-eigenpairs are the eigenpairs defined in  \cite{Qi05},
which satisfy
\be
\label{Heig}
\cA x^{m-1} = \lambda x^{[m-1]},
\ee
where $x^{[m-1]}= [x_{1}^{m-1}, x_{2}^{m-1}, \cdots, x_{n}^{m-1}]^T$.

\item  If $\cA \in \bR^{[m,n]}$, $m=m'$ and $\cB =\cI$  is the unit tensor,  mode-1 $\cB_R$-eigenvalues are the H-eigenvalues defined in  \cite{Qi05}.

\end{itemize}

\begin{5}
\label{remark1} 
{\rm 
Theoretical properties of mode-1 eigenvalues of tensors such as the Perron-Frobenius theory (\cite{CPZ08, FGH13, YY10, YY11}) for nonnegative tensors can be parallelly developed for mode-$k$ eigenvalues. However, as Horn and Johnson  indicated in \cite{HJ13}: ``One should not dismiss left eigenvectors as merely a parallel theoretical alternative to right eigenvectors. Each type of eigenvector can convey different information about a matrix,''  we believe that mode-1 through mode-$m$ eigenpairs can convey different information about a general tensor of order $m \geq 3$. 
}
\end{5}

In the rest of this section, we will obtain an upper bound for the number of equivalence classes of mode-$k$ eigenpairs.
As shown in Definition~\ref{def1}, Remark~\ref{remark4} and Remark~\ref{remark5}, the number of equivalence classes of mode-$k$ generalized eigenpairs for general tensors $\cA \in \bC^{[m,n]}$ and $\cB \in \bC^{[m',n]}$ is equivalent to the number of solutions to the following system of polynomials
\be \label{polyhyperplane}
T(\lambda,x)=
\begin{pmatrix}
(\cA^{(k)}x^{m-1})_1 - \lambda (\cB x^{m'-1})_1\\
\vdots \\
(\cA^{(k)}x^{m-1})_n - \lambda (\cB x^{m'-1})_n \\
a_1 x_1 + a_2x_2+ \cdots + a_n x_n+b
\end{pmatrix}
= 0,
\ee
where $\lambda$ and $x := (x_1, \cdots, x_n)^T$ are the unknowns, $a_1,\dots,a_n,b$ are random complex numbers. This 
motivates us to use Bernstein's theorem and its extensions in the field of solving polynomial systems (see \cite{Bernstein75, LW96}) to study the number of equivalence classes of eigenpairs. 

To initiate our discussion, we first introduce some commonly used notations and definitions. Let $P(x):=(p_1(x), \dots, p_n(x))^T$ be a polynomial system with $x:= (x_1, \dots, x_n)^T$. For $\alpha := (\alpha_1, \dots, \alpha_n)\in (\bZ_{\ge 0}^n)^T$, write $x^{\alpha} := x_1^{\alpha_1}\cdots x_n^{\alpha_n}$ and denote $|\alpha| = \alpha_1 + \dots + \alpha_n$. Then $P(x)$ can be denoted by
\be
\label{gpolsys}
P(x):= 
\begin{pmatrix}
p_1(x) :=\sum\limits_{\alpha \in S_1} c_{1,\alpha} x^{\alpha} \\
\vdots \\
p_n(x) :=\sum\limits_{\alpha \in S_n} c_{n,\alpha} x^{\alpha}
\end{pmatrix},
\ee
where $S_1, \dots, S_n$ are given finite subsets of $(\bZ_{\ge 0}^n)^T$ and $c_{i,\alpha} \in \bC^*:= \bC \backslash \{0\}$ are given coefficients of the corresponding monomials. Here for each $i = 1, \dots, n$, $S_i$ is called the support of $p_i(x)$ and its convex hull $R_i:=\mathrm{conv}(S_i)$ in $\bR^n$ is called the Newton polytope of $p_i(x)$. $(S_1, \dots, S_n)$ is called the support of $P(x)$. For nonnegative variables $\lambda_1, \dots, \lambda_n$, let $\lambda_1R_1 + \dots + \lambda_n R_n$ be the \emph{Minkowski} sum of $\lambda_1R_1, \dots, \lambda_n R_n$, i.e., 
$$
\lambda_1R_1 + \dots + \lambda_n R_n := \{\lambda_1 r_1 + \dots + \lambda_n r_n \, | \, r_i \in R_i, i = 1, \dots, n\}.
$$
 The $n$-dimensional volume of $\lambda_1 R_1 + \dots + \lambda_n R_n$, denoted by $\mathrm{Vol}_n(\lambda_1 R_1 + \dots + \lambda_n R_n)$, is a homogeneous polynomial function of degree $n$ in $\lambda_1,\dots,\lambda_n$ (See, for example, Proposition 4.9 of \cite{CLO05} for a proof). The coefficient of the monomial $\lambda_1\lambda_2\dots \lambda_n$ in $\mathrm{Vol}_n(\lambda_1 R_1 + \dots + \lambda_n R_n)$  is called the mixed volume of $R_1,\dots,R_n$, denoted by $\mathrm{MV}_n(R_1,\dots,R_n)$, or the mixed volume of the supports $S_1,\dots,S_n$, denoted by $\mathrm{MV}_n(S_1,\dots,S_n)$. Sometimes it is also called the mixed volume of $P(x)$ if no ambiguity exists. The following theorem relates the number of solutions of a polynomial system to its mixed volume.

\begin{2}
\label{theorem4}
\emph{(\textbf{Bernstein's Theorem}) \cite{Bernstein75} } The number of isolated zeros in $(\bC^*)^n$, counting multiplicities, of a polynomial system $P(x)=(p_1(x),\dots,p_n(x))^T$ with supports $S_1,\dots,S_n$ is bounded by the mixed volume $\mathrm{MV}_n(S_1,\dots,S_n)$. Moreover, for generic choices of the coefficients in $p_i$, the number of isolated zeros is exactly $\mathrm{MV}_n(S_1,\dots,S_n)$.
\end{2}

An unexpected limitation of Theorem~\ref{theorem4} is that it only counts the isolated zeros of a polynomial system in $(\bC^*)^n$ rather than $\bC^n$. To deal with this issue, Li and Wang gave the following theorem.

\begin{2}
\label{theorem5}
\emph{\cite{LW96}} The number of isolated zeros in $\bC^n$, counting multiplicities, of a polynomial system $P(x)=(p_1(x),\dots,p_n(x))^T$ with supports $S_1,\dots,S_n$ is bounded by the mixed volume $\mathrm{MV}_n(S_1 \cup \{0\},\dots,S_n \cup \{0\})$. 
\end{2}

The following lemma was given as Exercise 7 on page 338 of \cite{CLO05}.

\begin{1}
\label{lemma1}
Consider a polynomial system $P(x)=(p_1(x),\dots,p_n(x))^T$ with supports $S_1 = S_2 = \dots = S_n = S$. Then 
$$\mathrm{MV}_n(S,\dots,S) = n ! \mathrm{Vol}_n(\mathrm{conv}(S)).$$
\end{1}

An upper bound for the number of equivalence classes of mode-$k$ eigenpairs which generalizes results in \cite{CS13, OO13, Qi05} is given in the following theorem.

\begin{2}
\label{theorem2}
 Let $\cA \in \bC^{[m,n]}$ and $\cB \in  \bC^{[m',n]}$. Assume that  $ \cB x^{m'}$ is not identically zero as a function of $x$. Let $k$ be an integer such that $ 1 \leq k \leq m$. Assume that $\cA$ has finitely many equivalence classes of mode-$k$ $\cB$-eigenpairs over $\bC$. 
\begin{itemize}
\item[(a)] If $m = m'$, then the number of equivalence classes of mode-$k$ $\cB$-eigenpairs, counting multiplicities, is bounded by 
$$
n(m-1)^{n-1}.
$$
If $\cA$ and $\cB$ are generic tensors, then $\cA$ has exactly $n(m-1)^{n-1}$ equivalence classes of mode-$k$ $\cB$-eigenpairs.
\item[(b)] If $m \ne m'$, then the number of equivalence classes of mode-$k$ $\cB$-eigenpairs, counting multiplicities, is bounded by
$$
\frac{(m-1)^n-(m'-1)^n}{m-m'}.
$$
If $\cA$ and $\cB$ are generic tensors, then $\cA$ has exactly $((m-1)^n-(m'-1)^n)/(m-m')$ equivalence classes of mode-$k$ $\cB$-eigenpairs.
\end{itemize}
\end{2}
\bp 
Recall that the number of equivalence classes of mode-$k$ $\cB$-eigenpairs of $\cA$ is equal to the number of solutions of (\ref{polyhyperplane}).  For the random hyperplane $a_1x_1+\dots +a_nx_n +b = 0$ in (\ref{polyhyperplane}), without loss of generality, suppose that $a_n \ne 0$. Then $x_n$ can be solved as
\be
\label{xn}
x_n = c_1 x_1 + \dots + c_{n-1}x_{n-1} + d ,
\ee
where $c_i = -a_i/a_n$ for $i = 1,\dots, n-1$ and $d= -b/a_n$. Notice that the number of solutions of (\ref{polyhyperplane}) in $\bC^{n+1}$ is the same as the number of solutions in $\bC^n$ of the resulting system $T^*(\lambda,x_1,\dots,x_{n-1})$ by substituting (\ref{xn}) into the first $n$ equations of (\ref{polyhyperplane}). Denote the corresponding supports of $T^*$ by $S_1,\dots,S_n$. We claim that 
\be
\label{MV}
\mathrm{MV}_n(S_1 \cup \{0\},\dots,S_n\cup \{0\}) \le 
\left\{
  \begin{array}{ll}
    n(m-1)^{n-1}, & m=m' \\
   \dfrac{(m-1)^n - (m'-1)^n}{m - m'}, & m \ne m'
  \end{array}
\right.
\ee
Let $N$ denote the number of equivalence classes of mode-$k$ $\cB$-eigenpairs of $\cA$ over $\bC$. Then (\ref{MV}) implies that 
$$N\le  n(m-1)^{n-1}$$
for $m = m'$ and
$$ N \le \frac{(m-1)^n - (m'-1)^n}{m - m'}$$
for $m \ne m'$. When $\cA$ and $\cB$ are generic, equality holds in the two inequalities above by using Theorem~\ref{theorem4} and Theorem~\ref{theorem5}.

To prove (\ref{MV}), let $\bar{\cA} \in \bC^{[m,n]}$ and $\bar{\cB} \in \bC^{[m',n]}$ be generic tensors. Similar to (\ref{polyhyperplane}) the corresponding polynomial system to solve is
\be \label{randpolyhyperplane}
\bar{T}(\lambda,x)=
\begin{pmatrix}
(\bar{\cA}^{(k)}x^{m-1})_1 - \lambda (\bar{\cB} x^{m'-1})_1\\
\vdots \\
(\bar{\cA}^{(k)}x^{m-1})_n - \lambda (\bar{\cB} x^{m'-1})_n \\
a_1 x_1 + a_2x_2+ \cdots + a_n x_n+b
\end{pmatrix}
= 0.
\ee
Substituting (\ref{xn}) into the first $n$ equations of (\ref{randpolyhyperplane}) yields a new system $\bar{T}^*(\lambda,x_1,\dots,x_{n-1})$. The coefficients of the polynomials in the new system $\bar{T}^*$ are the sums of products of certain coefficients of the old system $\bar{T}$. For example, when $k=1$, the coefficient of $x_1^{m-1}$ in the first polynomial of the new system $\bar{T}^*$ is
\be \label{newcoef}
\sum\limits_{i=0}^{m-1} \left(\sum\limits_{\sigma \in I_i} a_{1,\sigma}\right)c_1^{i},
\ee
where $I_i$ is the set of all permutations of the set consisting of $i$ numbers of $n$ and $m-1-i$ numbers of $1$. Since $a_{1,\sigma}$'s and $c_1$ are all generic, the coefficient (\ref{newcoef}) is nonzero. Similarly, all other coefficients in the new system are nonzero. Now let $\bar{S}_1,\dots,\bar{S}_n$ be the corresponding supports of $\bar{T}^*$, we can assume that all monomials 
$$
\{x_1^{\alpha_1}x_2^{\alpha_2}\dots x_{n}^{\alpha_n} \Big{|} \alpha_i \in \bZ_{\ge 0},\, \alpha_1+\alpha_2+\dots+\alpha_n = m-1\}
$$ 
and 
$$
\{\lambda x_1^{\alpha_1}x_2^{\alpha_2}\dots x_n^{\alpha_n} \Big{|} \alpha_i \in \bZ_{\ge 0},\, \alpha_1+\alpha_2+\dots+\alpha_{n} = m'-1\}
$$ 
will appear in each of the first $n$ equations in (\ref{randpolyhyperplane}). Therefore, after substituting (\ref{xn}) into the first $n$ equations of (\ref{randpolyhyperplane}), all monomials 
$$
\{x_1^{\alpha_1}x_2^{\alpha_2}\dots x_{n-1}^{\alpha_{n-1}} \Big{|} \alpha_i \in \bZ_{\ge 0},\, \alpha_1+\alpha_2+\dots+\alpha_{n-1} \le m-1\}
$$ 
and 
$$
\{\lambda x_1^{\alpha_1}x_2^{\alpha_2}\dots x_{n-1}^{\alpha_{n-1}} \Big{|} \alpha_i \in \bZ_{\ge 0},\, \alpha_1+\alpha_2+\dots+\alpha_{n-1} \le m'-1\}
$$ 
will be contained in each equation of $\bar{T}^*$. This implies that $\bar{S}_1, \dots, \bar{S}_n$ are all equal to 
$$
\bar{S} := \{(0,\alpha) \,\big{|}\, \alpha \in (\bZ_{\ge 0}^{n-1})^T, |\alpha| \le m-1\} \cup \{(1,\alpha)\,\big{|} \, \alpha \in (\bZ_{\ge 0}^{n-1})^T, |\alpha| \le m'-1\}.
$$
 Notice that the convex hull of the set $\{\alpha \in (\bZ_{\ge 0}^{n-1})^T \, \big{|} \, |\alpha| \le m-1 \}$
is the $(n-1)$-simplex in $R^{n-1}$ with vertices $(0,0, \cdots,0), (m-1,0,\cdots,0), \cdots$, $(0, \cdots,0, m-1)$,  and the convex hull of the set $\{\alpha \in (\bZ_{\ge 0}^{n-1})^T \, \big{|} \, |\alpha| \le m'-1 \}$ is  the $(n-1)$-simplex in $R^{n-1}$ with vertices $(0,0, \cdots,0), (m'-1,0,\cdots,0), \cdots, (0, \cdots, 0, m'-1)$.
 Thus, their volumes are $(m-1)^{n-1}/(n-1)!$ and $(m'-1)^{n-1}/(n-1)!$ respectively (see, for example, Exercises 2 and 3 on page 307 of \cite{CLO05}).  Let $\bar{Q}$ be the convex hull of $\bar{S}$.
Then $\bar{Q}$ is the linear interpolation between the two aforementioned simplices with 
$$\bar{Q}= \{(t,\mathrm{conv}(\alpha)) \, \big{|} \, t \in[0,1], \alpha \in (\bZ_{\ge 0}^{n-1})^T, |\alpha| \le m-1+t(m'-m)\}.$$
Since for fixed $t$, $\mathrm{conv}(\{\alpha \in (\bZ_{\ge 0}^{n-1})^T \, \big{|} \, |\alpha| \le m-1+t(m'-m)\})$ is a simplex with volume $(m-1+t(m'-m))^{n-1}/(n-1)!$, we have
$$
\mathrm{Vol}_n(\bar{Q}) 
= \int_0^1 \frac{(m-1 + (m'-m)t)^{n-1}}{(n-1)!} dt
=
\left\{
\begin{array}{ll}
	\dfrac{(m-1)^n-(m'-1)^n}{(m-m')n!}, & m\ne m', \\
	\dfrac{n(m-1)^{n-1}}{n!}, & m=m'.  \\
\end{array}
\right.
$$
Therefore, by Lemma~\ref{lemma1},
$$
\mathrm{MV}_n(\bar{S_1},\dots,\bar{S}_n) = n! \mathrm{Vol}_n(\bar{Q}) = \left\{
\begin{array}{ll}
\dfrac{(m-1)^n-(m'-1)^n}{m-m'}, & m\ne m', \\
n(m-1)^{n-1}, & m=m'.  \\
\end{array}
\right.
$$

When tensors $\cA$ and $\cB$ in (\ref{polyhyperplane}) have some zero entries,  each support $S_i$ of $T^*$ must be a subset of the support $\bar{S}_i$ of $\bar{T}^*$, i.e., $S_i\subset \bar{S}_i$. In addition, since each polynomial of $\bar{T}^*$ contains a constant term (which arises when substituting (\ref{xn}) into the term $x_n^{m-1}$ in the first $n$ polynomials of $\bar{T}$), $0 \in \bar{S}_i$. Thus we have $S_i\cup \{0\} \subset \bar{S}_i$. Since mixed volume is monotonic \cite{Burago01},
$$
\mathrm{MV}_n(S_1\cup \{0\}, \dots, S_n\cup \{0\}) \le \mathrm{MV}_n(\bar{S_1},\dots,\bar{S}_n).
$$
This implies that (\ref{MV}) holds.

\ep

\begin{5} 
\label{remark2} {\rm
A few remarks about Theorem~\ref{theorem2}: 
\begin{itemize}
\item[(a)]  Theorem~\ref{theorem2} provides an upper bound for the number of equivalence classes of $\cB$-eigenpairs of tensor $\cA$ if $\cA$ has finitely many $\cB$-eigenpairs. The bound is tight when $\cA$ and $\cB$ are generic. In the literature, the numbers of eigenvalues as defined in (\ref{Heig}) and E-eigenpairs have been investigated.  

When $m=m'$ and $\cB$ is the $m$th order, $n$-dimensional identity tensor, Qi \cite[Theorem 1]{Qi05} and Chang, Qi, and Zhang \cite[Remarks 1]{CQZ13} proved that a tensor $\cA \in \bR^{[m,n]}$ has  $n(m-1)^{n-1}$ eigenvalues (\ref{Heig}).  In this case, Part (a) of Theorem~\ref{theorem2} gives  $n(m-1)^{n-1}$ as the upper bound for the number of equivalence classes of such eigenpairs.  

 When $m'=2$ and $\cB$ is the $n \times n$ identity matrix, Cartwright and Sturmfels \cite[Theorem 1.2]{CS13} proves that tensor $\cA$ has $((m-1)^n-1)/(m-2)$ equivalent classes of E-eigenpairs if it has finitely many E-eigenpairs.  Part (b) of our Theorem~\ref{theorem2} gives $((m-1)^n-1)/(m-2)$ as  the upper bound of the number of equivalence classes of E-eigenpairs if $\cA$ has finitely many E-eigenpairs.

\item[(b)] The upper bound given in Theorem~\ref{theorem2} can be highly useful in designing effective homotopy methods for computing mode-$k$ generalized eigenpairs. In fact, the homotopy method described in Algorithm~\ref{algorithm1} for the case $m=m'$ relies on the bound $n(m-1)^{n-1}$.   

\end{itemize}
}
\end{5}

\section{Computing complex tensor eigenpairs via homotopy methods}
\label{sec3}
\setcounter{equation}{0}

Consider $\cA \in \bC^{[m,n]}$ and $\cB \in \bC^{[m',n]}$. As discussed in Section~\ref{sec2}, the problem of computing mode-$k$ $\cB$-eigenpairs of $\cA$ in (\ref{tevalue}) or (\ref{evalue}) is equivalent to the problem of solving (\ref{polyhyperplane}), and if $m\ne m'$, normalize $(\lambda,x)$ to satisfy that $\cB x^{m'}=1$. Since (\ref{polyhyperplane}) is  a polynomial system, we consider to use a homotopy continuation method to numerically solve it. 

The basic idea of using homotopy continuation method to solve a general polynomial system $P(x)=(p_1(x),\dots,p_n(x))^T=0$ as defined in (\ref{gpolsys}) is to first deform $P(x)=0$ to another polynomial system $Q(x)=0$ that is easy to solve. Specifically, we construct a homotopy $H: \bC^n \times [0,1] \to \bC^n$ such that $H(x,0) = Q(x)$ and $H(x,1) = P(x)$. Then under certain conditions, the homotopy $H(x,t)=0$ has smooth solution paths parameterized by $t$ for $t \in [0,1)$ and all the isolated solutions of $P(x) = 0$ can be reached by tracing these paths.

A useful homotopy is the linear homotopy (see \cite{Li03,Morgan09, SW05, Wright85}):
\be
\label{linhom}
H(x,t) = (1-t)\gamma Q(x)+tP(x) = 0, \quad t \in [0,1],
\ee
where $\gamma$ is a generic nonzero complex number. It is very critical to choose a suitable $Q(x)$ such that the system $Q(x)=0$ is easy to solve and all isolated solutions of $P(x)=0$ can be found by tracing solution curves of $H(x,t)=0$.    

One choice of $Q(x)$ that always makes the linear homotopy (\ref{linhom}) work is the so-called total degree homotopy, in which the starting system $Q(x)=0$  has $deg=d_1\times d_2 \times \cdots \times d_n$ solutions (\cite{Li03,Morgan09,Wright85}), where $d_1,\dots,d_n$ are the degrees of polynomials $p_1(x), \dots, p_n(x)$ respectively.  $deg$ is called the total degree or  B\'{e}zout's number. By tracking $deg$ number of solution paths of (\ref{linhom}) we can find all the isolated solutions of $P(x)=0$. However, most  polynomial systems in applications usually have far fewer than $deg$ solutions. In this case, many of the $deg$  paths will diverge to infinity as $t \to 1$ resulting in huge wasteful computations. 

The polyhedral homotopy \cite{HS95} based on Bernstein's Theorem \cite{Bernstein75} makes significant progress in this sense. In this method, the number of paths that need to be traced is the mixed volume of a polynomial system, which generally provides a much tighter bound than B\'{e}zout's number for the number of isolated zeros of a polynomial system. Hence the new method reduces a significant amount of extraneous paths than the total degree homotopy in most occasions and thereby is much more efficient. However, the polyhedral homotopy includes two major stages: mixed volume computation and tracking paths. The computation of mixed volumes is a sophisticated procedure \cite{BHSW}. Moreover, mixed volume computation can be very expensive  for large polynomial systems. Thus if the mixed volume is far less than the B\'{e}zout's number and an appropriate linear homotopy can be constructed so that only mixed volume number of paths need to be traced, the system is better to be solved by using a linear homotopy instead of the polyhedral homotopy. 

To compute tensor eigenpairs, one can certainly use the polyhedral homotopy implemented in HOM4PS \cite{LLT08},  PHCpack \cite{Verschelde99}, PHoM \cite{GKKTFM04},  PSOLVE \cite{ZL13} (which is a MATLAB implementation of HOM4PS),  or the total degree homotopy implemented in Bertini \cite{BHSW}. However, using these methods to solve (\ref{polyhyperplane}) or (\ref{evalue}) directly does not take advantage of the special structures of a tensor eigenproblem. We will introduce two homotopy-type algorithms here that utilize such structures.

\subsection{A linear homotopy method when $m=m'$}

Theorem~\ref{theorem2} gives us that the mixed volume of (\ref{polyhyperplane}) when $m=m'$ is $n(m-1)^{n-1}$, which is far less than the B\'{e}zout's number, $m^n$. We consider constructing a linear homotopy in which only the mixed volume number of paths are traced.

 For a polynomial system $P(x)=(p_1(x),\dots,p_n(x))^T$ as defined in (\ref{gpolsys}), where $x = (x_1,\dots,x_n)$. Partition the variables $x_1,\dots,x_n$ into $k$ groups $y_1 = (x_1^{(1)},\dots,x_{l_1}^{(1)}), y_2 = (x_1^{(2)},\dots,x_{l_2}^{(2)}), \dots, y_k = (x_1^{(k)},\dots,x_{l_k}^{(k)})$ with $l_1 + \dots + l_k = n$. Let $d_{ij}$ be the degree of $p_i$ with respect to $y_j$ for $i = 1,\dots,n$ and $j=1,\dots,k$. Then the multihomogeneous B\'{e}zout's number of $P(x)$ with respect to $(y_1,\dots,y_k)$ is the coefficient of $\alpha_1^{l_1}\alpha_2^{l_2}\dots\alpha_k^{l_k}$ in the product
$$
\prod_{i=1}^n (d_{i1}\alpha_1 + \dots + d_{ik}\alpha_k).
$$

The following theorem will play an important role in constructing a proper linear homotopy.

\begin{2}
\label{theorem6}
\emph{\cite{SW05}} Let $Q(x)$ be a system of polynomials chosen to have the same multihomogeneous form as $P(x)$ with respect to certain partition of the variables $(x_1,\dots,x_n)$. Assume $Q(x)=0$ has exactly the multihomogeneous B\'{e}zout's number $\mathcal{N}$ of nonsingular solutions with respect to this partition. Then for almost all $\gamma \in \bC^*$, the homotopy 
$$
H(x,t)= (1-t) \gamma Q(x) + t P(x) = 0,
$$
has $\mathcal{N}$ nonsingular solution paths on $t \in [0,1)$ whose endpoints as $t \to 1$ include all the isolated solutions of $P(x)=0$.
\end{2}

For (\ref{polyhyperplane}), when $m = m'$ the following polynomial system
\be \label{qieig}
G(\lambda,x)=
\begin{pmatrix}
(\cA^{(k)}x^{m-1})_1 - \lambda (\cB x^{m-1})_1\\
\vdots \\
(\cA^{(k)}x^{m-1})_n - \lambda (\cB x^{m-1})_n \\
a_1 x_1 + a_2x_2+ \cdots + a_n x_n+b
\end{pmatrix}
= 0
\ee
needs to be solved, where $\lambda$ and $x := (x_1, \cdots, x_n)^T$ are the unknowns, $a_1,\dots,a_n,b$ are random complex numbers. Consider the starting system
\be \label{startQ}
Q(\lambda,x) = 
\begin{pmatrix}
(\lambda- \mu_1) (x_1^{m-1}-\beta_1)\\
(\lambda- \mu_2) (x_2^{m-1}-\beta_2)\\
\vdots \\
(\lambda- \mu_n) (x_n^{m-1}-\beta_n)\\
c_1x_1 + \dots c_n x_n + d
\end{pmatrix} = 0,
\ee
where $\mu_i, \beta_i, c_i$ for $i=1,\dots,n$ and $d$ are random nonzero complex numbers. 

\begin{2}
\label{theorem7}
Let $G(\lambda,x)$ and $Q(\lambda,x)$ be defined as (\ref{qieig}) and (\ref{startQ}) respectively. Then all the isolated zeros $(\lambda,x)$ in $\bC^{n+1}$ of $G(\lambda,x)$ can be found by using the homotopy
\be \label{mhomotopy}
H(\lambda,x,t)= (1-t) \gamma Q(\lambda,x) + t G(\lambda,x) = 0, \quad t \in [0,1]
\ee
for almost all $\gamma \in \bC^*$.
\end{2}
\bp
It is sufficient to verify that $Q(\lambda,x)$ satisfies all the assumptions of Theorem~\ref{theorem6}.
Partition the variables $(\lambda,x_1,\dots,x_n)$ into two groups: $(\lambda)$ and $(x_1,\dots,x_n)$, we can easily see that each of the first polynomial equations in (\ref{qieig}) and (\ref{startQ}) has degree 1 in $(\lambda)$ and degree $m-1$ in $(x_1,\dots,x_n)$, and the last equation in both (\ref{qieig}) and (\ref{startQ}) has degree 0 in $(\lambda)$ and degree $1$ in $(x_1,\dots,x_n)$. Hence (\ref{qieig}) and (\ref{startQ}) have the same multihomogeneous B\'{e}zout's number as the coefficient of $\alpha_1 \alpha_2^n$ in the product
$$
[1 \cdot \alpha_1 + (m-1)\alpha_2]^n(0\cdot \alpha_1 + 1 \cdot \alpha_2).
$$
It can be easily computed that this coefficient is equal to 
$$
\begin{pmatrix}
n \\
1
\end{pmatrix}
(m-1)^{n-1}
=n(m-1)^{n-1}.
$$
Hence (\ref{qieig}) and (\ref{startQ}) have the same multihomogeneous B\'{e}zout's number $n(m-1)^{n-1}$ with respect to the partition $(\lambda)$ and $(x_1,\dots,x_n)$.

We now show that $Q(\lambda,x)$ in (\ref{startQ}) has exactly $n(m-1)^{n-1}$ zeros. Notice that if $\lambda$ is equal to none of $\mu_1,\dots,\mu_n$, then we end up with a system of $n+1$ equations and $n$ unknowns, which has no solutions. Thus $\lambda$ must be equal to one of $\mu_1,\dots, \mu_n$. Without loss of generality, assume that $\lambda=\mu_1$. Then $x_1,\dots,x_n$ can be determined by
\bqn
x_i^{m-1}-\beta_i &=& 0, \quad i = 2,\dots,n\\
c_1x_1 + \dots c_n x_n + d &=& 0
\eqn
Obviously, each $x_i$ for $i=2,\dots,n$ can be chosen as one of the $(m-1)$-th root of $\beta_i$ and $x_1$ will be solved by substituting the chosen $x_2,\dots,x_n$ into the last hyperplane equation. So there are $(m-1)^{n-1}$ solutions corresponding to $\lambda=\mu_1$. This argument holds for $\lambda$ being any $\mu_i$. Therefore, there are totally $n(m-1)^{n-1}$ solutions.

It remains to prove that each solution of $Q(\lambda,x)=0$ in (\ref{startQ}) is nonsingular. As discussed above, any solution $(\lambda^*,x^*)$ of (\ref{startQ}) satisfies 
\bea
\lambda^* &=& \mu_i, \nonumber\\
(x_{j}^*)^{m-1} - \beta_{j} &=& 0, \ \ j=1, \cdots, i-1, i+1, \cdots, n,        \label{Qsol}\\
c_1x_1^* + \dots + c_nx_n^* + d &=& 0. \nonumber
\eea
Let $DQ(\lambda,x)$ be the Jacobian of $Q(\lambda,x)$ with respect to $(\lambda,x)$. It is sufficient to show that $DQ(\lambda^*,x^*)$ is nonsingular. Denote 
$$A_j(\lambda,x) := x_j^{m-1} - \beta_j, \quad B_j(\lambda,x) := (\lambda-\mu_j)(m-1)x_j^{m-2}$$
for $j = 1,\dots,n$. Then 
$$
DQ(\lambda,x) = 
\begin{pmatrix}
A_1 & B_1 &  &  &  &  & & \\
\vdots & & \ddots & &  & & & \\
A_{i-1} &  & & B_{i-1} & & & & \\
A_i & & & & B_i & & & \\
A_{i+1} & & & & &  B_{i+1} & & \\
\vdots & & & & &  & \ddots & \\
A_n & & & & &  &  & B_n\\
0 & c_1 & \dots & c_{i-1} & c_i & c_{i+1} & \dots & c_n
\end{pmatrix}.
$$
Note that 
$A_j(\lambda^*,x^*) = (x_j^*)^{m-1} - \beta_j= 0, \quad j \ne i$
and 
$B_i(\lambda^*,x^*) = (\lambda^*-\mu_i)(m-1)(x_i^*)^{m-2}= (\mu_i-\mu_i)(m-1)(x_i^*)^{m-2}=0$
by (\ref{Qsol}). For simplicity, write $A_j^*:=A_j(\lambda^*,x^*) $ and $B_j^*:=B_j(\lambda^*,x^*)$. Then
$$
DQ(\lambda^*,x^*) = 
\begin{pmatrix}
0 & B_1^* &  &  &  &  & & \\
\vdots & & \ddots & &  & & & \\
0 &  & & B_{i-1}^* & & & & \\
A_i^* & & & & 0& & & \\
0 & & & & &  B_{i+1}^* & & \\
\vdots & & & & &  & \ddots & \\
0 & & & & &  &  & B_n^*\\
0 & c_1 & \dots & c_{i-1} & c_i & c_{i+1} & \dots & c_n
\end{pmatrix}.
$$
Then
$$
\det(DQ(\lambda^*,x^*)) = (-1)^{i+1}A_i^*(-1)^{n+i}c_i\prod_{j \ne i} B_j^* \ne 0
$$
by (\ref{Qsol}).
\ep

Theorem~\ref{theorem7} suggests us that (\ref{mhomotopy}) can be used to solve (\ref{polyhyperplane}) in the case of $m = m'$. For simplicity, write $u:= (\lambda,x)$. In order to improve numerical stability, we apply the transformation $s = \ln t$ to (\ref{mhomotopy}) (a strategy first suggested in \cite{GKKTFM04}) and obtain the new homotopy as
\be \label{mhomotopy-s}
\bar{H}(u,s)= (1-e^s) \gamma Q(u) + e^s G(u) = 0, \quad s\in [-\infty,0]
\ee
where $Q(u)=Q(\lambda,x)$ and $G(u)=G(\lambda,x)$ are defined in (\ref{startQ}) and (\ref{qieig}) respectively. 

We now introduce our linear homotopy method for computing mode-$k$ generalized eigenpairs when $m=m'$. 

\begin{70}
\label{algorithm1}
{\rm (Compute  complex mode-$k$ $\cB$-eigenpairs of $\cA$, where $\cA, \cB \in \bC^{[m,n]}$.)} \\

{\bf Step 1.} Compute all solutions of $Q(u)$ as defined in (\ref{startQ}). 

{\bf Step 2.} Path following: Follow the paths from $s = -\infty$ to $s=0$. In reality, we certainly cannot start from $s= -\infty$. In this case, one can choose a very negative $s_0$ and obtain a starting point by using Newton's iterations:
$$
w^{(k+1)} = w^{(k)} - [\bar{H}_u(w^{(k)},s_0)]^{-1}\bar{H}(w^{(k)},s_0), \quad k = 0,1, \dots 
$$
until $\|\bar{H}(w^{(N)},s_0)\|$ is very small for some $N$. Here $w^{(0)}$ is a solution of $Q(u)=0$. Let $u_0 := u(s_0)$ and take $u_0 = w^{(N)}$. Then path following can be triggered.

Path following is  done  using the 
prediction-correction method. Let $(u_k, s_k):= (u(s_k), s_k)$, to find the next point on the path $\bar{H}(u,s) = 0$, we employ the following strategy:
\begin{itemize}
\item Prediction Step: Compute the tangent vector $\dfrac{du}{ds}$ to $\bar{H}(u,s)=0$ at $s_k$ by solving the linear system
\begin{equation*}
\bar{H}_u(u_k,s_k)\dfrac{du}{ds} = -\bar{H}_s(u_k,s_k)
\end{equation*}
for $\dfrac{du}{ds}$. Then compute the approximation 
$\tilde{u}$ to $u_{k+1}$ by
$$
\tilde{u} = u_k + \Delta s \frac{du}{ds}, \quad s_{k+1} = s_k+ \Delta s,
$$
where $\Delta s$ is a stepsize.
\item Correction Step: Use Newton's iterations. Initialize $v_0 = \tilde{u}$. For $i = 0, 1, 2, \dots$, compute
$$
v_{i+1} = v_i - [\bar{H}_u(v_i,s_{k+1})]^{-1}\bar{H}(v_i,s_{k+1})
$$
until $\|H(v_{J},s_{J})\|$ is very small. Then let $u_{k+1} = v_{J}$.
\end{itemize}

{\bf Step 3.} End game. During Step 3  when $s_N$ is very close to 0, the corresponding $u_N$ should be very close to a zero $u^*$ of $G(u)=G(\lambda,x)$. So Newton's iterations
$$
u^{(k+1)} = u^{(k)} - [DG(u^{(k)})]^{-1}G(u^{(k)}), \quad k = 0,1, \dots 
$$
will be employed to refine our final approximation $\tilde{u}$ to $u^*$. If $DG(u^*)$ is nonsingular, then $\tilde{u}$ will be a very good approximation of $u^*$ with multiplicity 1. If $DG(u^*)$ is singular,  $\tilde{u}$ is either an isolated singular zero of $G(u)$ with some multiplicity $l>1$ or  in a positive dimensional solution component of $G(u)=0$. We use a strategy provided in Chapter VIII of \cite{Li03} (see also \cite{SW05}) to verify whether $\tilde{u}$ is  an isolated zero with multiplicity $l>1$ or in a positive dimensional solution component of $G(u)=0$.

{\bf Step 4.} For each solution $u = (\lambda,x)$ obtained in Step 3, normalize $x$ in the following fashion to get a new eigenvector
\be
\label{normalization}
y=\frac{x}{x_{i_0}}
\ee
can be obtained, where $ i_0 := \operatorname*{arg\,max}_{1\le i \le n} |x_i|$. Hence $(\lambda,y)$ is an eigenpair.
\end{70}

\begin{5} 
\label{remark6} {\rm
A few remarks about Algorithm~\ref{algorithm1}:
\begin{itemize}
\item[(a)] As defined in (\ref{tevalue}) and Remark~\ref{remark4}, if $(\lambda,x)$ is an eigenpair, $(\lambda,tx)$ for $t \ne 0$ is also an eigenpair. Therefore, Step 4 is well defined in this sense.

\item[(b)] Notice that if $x$ is a real eigenvector associated with a real eigenvalue $\lambda$, $tx$ for any $t \in \bC \backslash \{0\}$ will be a complex eigenvector associated with the same eigenvalue $\lambda$. If in any case a complex eigenvector like $tx$ is obtained in Step 3 of Algorithm~\ref{algorithm1}, applying (\ref{normalization}) to $tx$ will give us a new real eigenvector. In this sense, Step 4 is very helpful for us to detect real eigenpairs.

\item[(c)]  According to Theorem \ref{theorem2}, if $\cA$ has finitely many equivalence classes of $\cB$-eigenpairs, then the  number of equivalence classes of $\cB$-eigenpairs,
counting multiplicities, is bounded  by $n(m-1)^{n-1}$. Moreover, this bound is attained when $\cA$ and $\cB$ are generic. This result implies that the optimal number of paths to follow in a homotopy method for solving the system (\ref{qieig}) is $n(m-1)^{n-1}$. Our starting system (\ref{startQ}) has exactly $n(m-1)^{n-1}$ nonsingular solutions. In this sense, 
Algorithm \ref{algorithm1} follows the optimal number of paths for solving the system (\ref{qieig}). 

\end{itemize}
}
\end{5}

\subsection{A polyhedral homotopy method when $m\ne m'$}

 To compute  mode-$k$ generalized tensor eigenpairs when $m \ne m'$, we use the equivalence class structure of the eigenproblem as described in Remark~\ref{remark4}.  We  first solve (\ref{polyhyperplane}) to find a representative $(\lambda, x)$ from each equivalence class  and then find all $m' $ eigenpairs from each equivalence class by simply using $\lambda'=t^{m-m'}\lambda, x'=tx$, where $t$ is a root of $t^{m'}=1$. We  use a polyhedral homotopy method to solve the system (\ref{polyhyperplane}). In our implementation,  the polyhedral homotopy method is  PSOLVE (\cite{ZL13}), with some modifications, as described in Subsection 3.3. The modifications are  
based on Strategies 2 and 3 introduced in the next subsection.

One may think of solving (\ref{evalue}) directly to get $m' $ eigenpairs from each equivalence class. However, this alternative method would have to follow $m'$ times as many paths as  the approach we described in the previous paragraph and therefore it would need much more computation.

Now we present our algorithm for computing mode-$k$ generalized eigenpairs when $m \ne m'$. 

\begin{70}
\label{algorithm2}
{\rm (Compute  complex mode-$k$ $\cB$-eigenpairs of $\cA$, where $\cA\in \bC^{[m,n]}, \cB \in \bC^{[m',n]}$ with $m \ne m'$.)} \\ 

{\bf Step 1.} Using modified PSOLVE to get all solutions $(\lambda,x)$ of (\ref{polyhyperplane}).

{\bf Step 2.} For each $(\lambda,x)$ obtained in Step 1, if $\cB x^{m'} \ne 0$, normalize it to get an eigenpair $(\lambda^*,x^*)$ by
$$
\lambda^* = \frac{\lambda}{(\cB x^{m'})^{(m-m')/m'}}, \quad x^* = \frac{x}{(\cB x^{m'})^{1/m'}}
$$
to satisfy (\ref{evalue}).

{\bf Step 3.}  Compute $m'$ equivalent eigenpairs $(\lambda',x')$ of $(\lambda^*,x^*)$ by $\lambda' = t^{m-m'} \lambda^*$ and $x' = tx^*$ with $t$ being a root of $t^{m'}=1$.

\end{70}

\subsection{Implementation tips}

By using random complex numbers in the formulation of homotopy, theoretically, with probability one the solution paths do not cross or go to infinity in the middle.  In practice, however, two paths may become very close to each other and the magnitude of some components of a solution curve may become very large during the procedure of path tracking. This causes various numerical difficulties such as missing solutions or losing efficiency or stability.  In our implementation of Algorithms \ref{algorithm1} and \ref{algorithm2}, we use the following strategies to address these issues.  We focus our discussion on Algorithms \ref{algorithm1}. Step 1 of Algorithm \ref{algorithm2} is done similarly. 

When tracing two paths that are sufficiently close, it is possible for the path tracing algorithm to jump from one path to the other path and thus result in the missing of zeros. To minimize the chance for curve jumping and keep the efficiency, our {\it First Strategy} is:  {\it  The stepsize $\Delta s$ in Step 2 of Algorithm \ref{algorithm1} is chosen adaptively}. Initially, set $s_0 = -20(n+1)$, where $n+1$ is the number of unknown variables $\lambda,x_1,\dots,x_n$ in (\ref{polyhyperplane}) or (\ref{qieig}), and $\Delta s = -s_0/3$. Similar to \cite{LLT08}, if more than 3 steps of Newton iterations were required to converge within the desired accuracy, $\Delta s$ is halved and the shorter step is attempted. On the other hand, if several (the default being 2) consecutive steps were not cut, $\Delta s$ is doubled, up to a prescribed maximum value (the default being $-s_k/3$). 

Although this adaptive approach can often significantly reduce the possibility of curve jumping, it can still occur in some cases. Our {\it Second Strategy} is: {\it To check if there is curve jumping, we store all the found solutions in a binary search tree. Each time when a new solution is found, we can quickly find (with time complexity $O(\log N)$, where $N$ is the number of solutions) whether there is any existing solution that is numerically identical to the new solution, that is, the difference between them is less than a threshold (the default being $10^{-6}$). If two numerically identical solutions are detected and the condition numbers of their Jacobian matrices are greater than a  threshold (the default being $10^{10}$), we consider that curve jumping likely has occurred.  We then retrace the two associated curves with more restrictively chosen parameters in the projective space, as described in the paragraph after the next one.}

When the magnitude of some components of certain solution curves become very large in the middle, tracing these paths may fail due to numerical instability. This issue can be largely resolved by following paths in the projective space (see, for example, \cite{SW05}). However, empirically it is more time consuming to follow all paths in the projective space than in the complex space. In our implementation of Algorithms \ref{algorithm1} and \ref{algorithm2}, our {\it Third Strategy} is: {\it To retrace solution curves in the projective space only for those paths that are detected to have very large solution components.}

To trace a path in the projective space, we  first homogenize each polynomial equation of (\ref{mhomotopy-s}) in the variables $\lambda,x_1,\dots,x_n$ to get the homotopy
\be \label{mhomotopy-p}
\footnotesize{
\hat{H}(\lambda,\hat{x},s) = (1-e^s) \gamma
\begin{pmatrix}
	(\lambda- \mu_1 x_0) (x_1^{m-1}-\beta_1 x_0^{m-1})\\
	\vdots \\
	(\lambda- \mu_n x_0) (x_n^{m-1}-\beta_n x_0^{m-1})\\
	c_1x_1 + \dots c_n x_n + d x_0
\end{pmatrix} + e^s
\begin{pmatrix}
	x_0(\cA^{(k)}x^{m-1})_1 - \lambda (\cB x^{m-1})_1\\
	\vdots \\
	x_0 (\cA^{(k)}x^{m-1})_n - \lambda (\cB x^{m-1})_n \\
	a_1 x_1 + a_2x_2+ \cdots + a_n x_n+b x_0
\end{pmatrix} = 0,}
\ee
where $\hat{x} = (x_0,x_1,\dots,x_n)^T$, and then follow  the solution curve of (\ref{mhomotopy-p}) in the projective space. Notice that in (\ref{mhomotopy-p}) if $(\lambda,\hat{x})$ is a solution, so is $(\alpha\lambda,\alpha\hat{x})$ for $\alpha \in \bC \backslash\{0\}$. Thus along the path we can always scale $(\lambda,\hat{x})$ to keep each component's magnitude in a suitable finite range.

To the best of our knowledge, Strategies 2 and 3 have not been used in other implementations of homotopy  methods, although some packages may trace all curves in the projective space.

\section{Computing real tensor eigenpairs via homotopy methods}
\label{sec4}
\setcounter{equation}{0}

In some applications, tensor $\cA$ is real and only real  eigenpairs (or real eigenvalues) of $\cA$ are of interest (\cite{CDN14, Qi05}).  In this situation, only real zeros of the polynomial system (\ref{evalue}) or (\ref{qieig}) are needed. It is worth noting that there is currently no effective method to find all real zeros for a polynomial system directly. One may think to use a homotopy continuation method to trace only real zeros from the start system to the target system. However,  this approach generally does not guarantee a real zero at the end, because the homotopy methods inherently  have to trace paths in the complex space in order to avoid the discriminant locus. 

For a real tensor $\cA$,  a real eigenvalue may have complex eigenvectors. Sometimes identifying which eigenvalues are real is the only concern. In this case, we can first compute complex zeros $(\lambda,x)$ of (\ref{qieig}) by Algorithm~\ref{algorithm1} or (\ref{evalue}) by Algorithm~\ref{algorithm2}, then identify the real eigenvalues by checking the size of the imaginary parts of $\lambda$'s. Specifically,  let $(\lambda^*, x^*)$ be a computed eigenpair. If
 $$
|\Im(\lambda^*)| < \delta_0,
$$ 
 where $\delta_0$ is a threshold for the imaginary part (the default value  $\delta_0=10^{-8}$), then we take $\Re(\lambda^*)$ as a real eigenvalue.  

Note that if $\frac{m'}{m-m'} $ is a nonzero integer multiple of 4 (for example, $m=5, m'=4$ or $m=10, m'=8$,) and $\cA$ has an eigenpair $(\lambda^*, x^*)$ with a purely imaginary eigenvalue $\lambda^*=bi$, where $b \in \bR$, then one can easily show that   $(b, (-i)^{1/(m-m')} x^*)$ and  $(-b, i^{1/(m-m')} x^*)$ are eigenpairs with real 
eigenvalues. Therefore, when $\frac{m'}{m-m'} $ is a nonzero integer multiple of 4, if $(\lambda^*, x^*)$ is an eigenpair found by Algorithm~\ref{algorithm2} such that 
$$
|\Re(\lambda^*)| < \delta_0,
$$ 
then we take $ \Im(\lambda^*)$ and $ -\Im(\lambda^*)$  as  real eigenvalues, with corresponding eigenvectors  $(-i)^{1/(m-m')} x^*$ and 
$i^{1/(m-m')} x^*$.

When looking for real tensor eigenpairs (i.e., both eigenvalues and eigenvectors being real), the situation becomes more complicated. We use a two-step procedure.  The first step is to get complex zeros $(\lambda,x)$ of (\ref{qieig}) by Algorithm~\ref{algorithm1} or (\ref{evalue}) by Algorithm~\ref{algorithm2}. The second step is to extract all real eigenpairs $(\lambda, x)$ from the complex zeros obtained in the first step. 

To facilitate the discussion, the following notation is introduced. For a vector $a=(a_1,\dots,a_n)^T \in \bC^n$, let
$$
\Im(a) = (\Im(a_1), \dots, \Im(a_n))^T, \quad \Re(a) = (\Re(a_1), \dots, \Re(a_n))^T.
$$

Suppose that $(\lambda^*, x^*)$ is an  eigenpair found in the first step. Consider two cases: (i)  $(\lambda^*, x^*)$ is an isolated eigenpair; (ii) $(\lambda^*, x^*)$ is an  eigenpair  contained in a positive dimensional solution component of system (\ref{qieig}) or (\ref{evalue}).

When $(\lambda^*, x^*)$ is an isolated eigenpair, take $(\Re(\lambda^*), \Re(x^*))$ as a real eigenpair if
 $$
\|\Im(\lambda^*, x^*)\|_2 < \delta_0.
$$ 

If $(\lambda^*, x^*) $  is an eigenpair in a positive dimensional solution component of system (\ref{qieig}) or (\ref{evalue}),  in general real eigenvectors are not guaranteed to be found by Algorithm~\ref{algorithm1} or Algorithm~\ref{algorithm2} even if the corresponding eigenvalue 
$\lambda^*$ is real. In this case, we will  
construct the following Newton homotopy \cite{AG90} 
\be
\label{Newtonhomotopy}
H(\lambda,x,t):= P(\lambda,x)-(1-t)P(\lambda^*,\Re(x^*)), \quad t \in [0,1]
\ee
to follow curves in $(\lambda,x)\in \bR^{n+1}$ in order to get a real eigenpair. Notice that when following curves in the complex space it is proved in \cite{Li03} that the solution curves of (\ref{Newtonhomotopy}) can be parameterized by $t$, but the solution curves of (\ref{Newtonhomotopy}) is not necessarily to be a function of $t$ when restricted in the real space. So a different method to follow curves is needed. In this case parameterizing the solution curves by the arc length $s$ is suggested in \cite{Li15}.  Instead of following paths using (\ref{Newtonhomotopy}), we will use the homotopy 
\be
\label{paraNewtonhomo}
H((\lambda (s),x(s),t(s))=0.
\ee 
For a description of the Newton homotopy method, we refer to \cite{Li15}.

An interesting phenomenon we have observed in our experiments is that in some  cases, the real eigenpairs can be obtained more straightforwardly from the complex eigenpairs found from Algorithm~\ref{algorithm1} or Algorithm~\ref{algorithm2}:  If $(\lambda^*,x^*)$ is in a positive dimensional solution component of (\ref{evalue}) and $\lambda^*\in \bR$, then $(\lambda^*,\Re(x^*)/ (\cB \Re(x^*))^{1/m'} )$ and $(\lambda^*,\Im(x^*)/ (\cB \Im(x^*))^{1/m'})$ can be mode-$k$ $\cB_R$ eigenpairs of $\cA$. 
This gives us a heuristic approach to find real eigenpairs for eigenpairs belong to  positive dimensional  components. We remark that this approach works well for all the examples (e.g., Example 4.8, 4.11, 4.13, 4.14) in \cite{CDN14} when a real Z-eigenvalue has infinitely many real Z-eigenvectors.   The following Proposition gives a justification for the approach in special cases.

\begin{4}
\label{lemma2}
Let $\cA \in \bR^{[m,n]}$ and $\cB \in  \bR^{[m',n]}$. Let $k$ be an integer such that $1 \leq k \leq m$. Let $\lambda\in \bR$ be a real mode-$k$ $\cB$ eigenvalue of $\cA$. If $U:=\{x \in \bC^n \, | \, \cA^{(k)} x^{m-1} = \lambda \cB x^{m'-1}\}$ contains a complex linear subspace $V$ of $\bC^n$ such that $y \in V$ implies $\bar{y} \in V$, then for any $x = \xi + i \eta \in V$ such that $\xi, \eta \in \bR^n$ and $\xi \ne 0, \eta \ne 0$, 
\begin{itemize}
\item When $m=m'$, $\xi$ and $\eta$ are  both real mode-$k$ $\cB$ eigenvectors of $\cA$ associated with $\lambda$. 
\item  When $m \ne m'$, $\cB \xi^{m'} \ne 0$ and $\cB \eta^{m'} \ne 0$, the normalized vectors
$$
v:=\frac{\xi}{(B\xi^{m'})^{1/m'}}, \quad w:=\frac{\eta}{(B\eta^{m'})^{1/m'}}
$$
 are  real mode-$k$ $\cB$ eigenvectors of $\cA$ associated with $\lambda$.

\end{itemize}

\end{4}
\bp
Let $x \in V$. Then $\bar{x} \in V$. Since $V$ is a linear subspace, $\xi=(x+\bar{x})/2$ and $\eta=(x-\bar{x})/(2i)$ are also in $V$. Thus, when $m = m'$,  $\xi$ and $\eta$ are  both real mode-$k$ $\cB$ eigenvectors of $\cA$ associated with $\lambda$. If $m \ne m'$, we have 
\bqn
\cB v^{m'} &=& \sum_{i_1, \cdots, i_{m'} =1}^{n} B_{i_1 i_2 \cdots i_{m'}}v_{i_1} v_{i_2} \cdots v_{i_{m'}}\\
&=& \sum_{i_1, \cdots, i_{m'} =1}^{n} B_{i_1 i_2 \cdots i_{m'}} \frac{\xi_{i_1}}{(B\xi^{m'})^{1/m'}} \frac{\xi_{i_2}}{(B\xi^{m'})^{1/m'}} \cdots \frac{\xi_{i_{m'}}}{(B\xi^{m'})^{1/m'}} \\
&=& \frac{\sum_{i_1, \cdots, i_{m'} =1}^{n} B_{i_1 i_2 \cdots i_{m'}}\xi_{i_1} \xi_{i_2} \cdots \xi_{i_{m'}}}{B\xi^{m'}} \\
&=& \frac{B\xi^{m'}}{B\xi^{m'}} = 1.
\eqn
This implies that $v$ is a real mode-$k$ $\cB$ eigenvector of  $\cA$ associated with $\lambda$. 
Similarly, $Bw^{m'}=1$ can also be verified. Therefore, $w$ is also a real mode-$k$ $\cB$ eigenvector of  $\cA$ associated with $\lambda$.
\ep

\begin{5}
\label{remark3} 
{\rm A natural question is when $U$ defined in Proposition \ref{lemma2} contains a linear subspace $V$. Consider the case when $\cA$ is a symmetric tensor with a low rank decomposition.  For simplicity, consider $m=3$. For vectors $a,b,c \in \bC^n$, define the outer product tensor $a \circ b \circ c = (a_i b_j c_k)$. Suppose that a symmetric tensor $\cA \in \bR^{[m,n]}$  can be decomposed as 
$$
\cA=   y_1 \circ y_1 \circ y_1 + y_2 \circ y_2 \circ y_2 + \cdots +  y_r \circ y_r \circ y_r,
$$
where $r<n$, $y_k \in \bR^n, k=1,2, \cdots, r$.   

Let $W={\rm span}(y_1, y_2, \cdots, y_r)$  and let $V$ be the orthogonal complement of $W$.  Then $V$ is a linear subspace of $\bC^n$. Clearly, $x \in V$ implies $\bar{x} \in V$. Moreover, for any $
x \in V \backslash \{0\}$,
$$
Ax^{m-1}=0.
$$
Thus $x$ is an eigenvector of $\cA$ corresponding to the eigenvalue $0$. In this case,   $\Re(x)/\|\Re(x)\|$ and $\Im(x)/\|\Im(x)\|$ are real Z-eigenvectors of $\cA$ corresponding to the real Z-eigenvalue $0$. Hence, the set $V$ is a desired linear subspace of $\bC^n$ contained in $U$.

}
\end{5}

Finally we present an algorithm for computing  real eigenpairs that combines the heuristic approach and the Newton homotopy method.

\begin{70}
\label{algorithm4}
{\rm (Compute real mode-$k$ $\cB$-eigenpairs of $\cA$, where $\cA\in \bR^{[m,n]}, \cB \in \bR^{[m',n]}$)}\\

{\bf Step 1.} Compute all complex eigenpairs  using Algorithm~\ref{algorithm1} or Algorithm~\ref{algorithm2}.  Let $K$ denote the set of 
  found eigenpairs $(\lambda, x)$ such that $|\Im(\lambda)| < \delta_0$.

{\bf Step 2.} For each eigenpair $(\lambda^*,x^*) \in K$:  if $(\lambda^*,x^*)$ is in a positive dimensional solution component of (\ref{qieig}) or (\ref{evalue}), go to Step 3. Otherwise, $(\lambda^*, x^*)$ is an isolated eigenpair. If $\|\Im( x^*)\|_2 < \delta_0$, then take $(\Re(\lambda^*), \Re(x^*))$ as a real eigenpair and stop.

{\bf Step 3.} Set $\lambda = \Re(\lambda^*)$. If $ m = m'$, set $v := \Re(x^*)$ (if $ \Re(x^*) \ne 0$)  and $w:=\Im(x^*)$   (if $ \Im(x^*) \ne 0$); otherwise,  set
$$ 
v:=\frac{\Re(x^*)}{(\cB\Re(x^*)^{m'})^{1/m'}}  \quad ({\rm if} \ \cB\Re(x^*)^{m'} \ne 0),
$$
and
$$   
w:=\frac{\Im(x^*)}{(\cB\Im(x^*)^{m'})^{1/m'}}  \quad ({\rm if} \ \cB\Im(x^*)^{m'} \ne 0).
$$
If $(\lambda, v)$ or $(\lambda,w)$ is a mode-$k$ $\cB$-eigenpair of $\cA$, then we have obtained a real eigenpair and stop. Otherwise, goto Step 4.   

{\bf Step 4.} Starting from $(\lambda^*,x^*)$, use the Newton homotopy method to follow curves of (\ref{paraNewtonhomo}) to find a real eigenpair.

\end{70}

\section{Numerical results}
\label{sec5}
\setcounter{equation}{0}

Based on the algorithms introduced in Section 3 and Section 4, a MATLAB  package \texttt{TenEig} has been developed.  The current version is \texttt{TenEig} 1.1. The numerical results reported in this paper were obtained using this version. The  package  can be downloaded from \\

\url{http://www.math.msu.edu/~chenlipi/TenEig.html}      \\

Consider the tensors $\cA \in \bC^{[m,n]}$ and $\cB \in \bC^{[m',n]}$. In the \texttt{TenEig} package, function \texttt{teig} can be used to compute the general mode-$k$ $\cB$  eigenvalues and eigenvectors of a tensor $\cA$ for $m = m'$. The  input of this function is: tensor $\cA$ or the polynomial form (if $\cA$ is symmetric) $\cA x^m$, tensor $\cB$ (the default is the identity tensor), mode $k$ (the default value is 1), and the output is: mode-$k$ $\cB$ eigenvalues and eigenvectors of $\cA$. By default, \texttt{teig} finds eigenvalues and eigenvectors for the eigenproblem (\ref{Heig}). 

The function \texttt{teneig} computes the general mode-$k$ $\cB$  eigenvalues and eigenvectors of a tensor $\cA$ for $m \ne m'$. The  input of this function is: tensor $\cA$ or the polynomial form (if $\cA$ is symmetric) $\cA x^m$, tensor $\cB$, mode $k$ (the default value is 1), and the output is: mode-$k$ $\cB$ eigenvalues and eigenvectors of $\cA$. If $\cB$ is chosen as the identify matrix, the \texttt{teneig} computes the E-eigenvalues and E-eigenvectors of $\cA$ as defined in Qi \cite{Qi05}. 
       
Since E-eigenpairs of a tensor are frequently needed, our package includes a separate function  \texttt{eeig}, which only computes E-eigenpairs of a tensor. 

The package also  includes two functions  \texttt{heig} and \texttt{zeig} to compute real eigenpairs of a tensor: The first one computes H-eigenpairs and the second one computes Z-eigenpairs.

In the next two subsections, numerical results are reported to illustrate the effectiveness and efficiency of our methods for computing tensor eigenpairs. All the numerical experiments were done on a Thinkpad T400 Laptop with an Intel(R) dual core CPU at 2.80GHz  and 2GB of RAM, running on a Windows 7 operating system. The package \texttt{TenEig} was run using  MATLAB 2013a. In our examples, we used \texttt{teig} or \texttt{teneig} to compute generalized eigenpairs, \texttt{teig} to compute eigenpairs (\ref{Heig}), \texttt{eeig} to compute E-eigenvalues, \texttt{heig} to compute (real) H-eigenpairs, and \texttt{zeig} to compute (real) Z-eigenpairs, respectively.

\subsection{Examples for computing complex eigenpairs}

In this subsection, some numerical examples illustrating the performance of \texttt{TenEig} for  computing complex tensor $\cB$-eigenpairs are provided.

A numerical solver \texttt{NSolve}, based on the Gr\"{o}bner basis, for solving systems of algebraic equations is provided by Mathematica. We will compare the performance of \texttt{TenEig} and \texttt{NSolve} in computing complex tensor $\cB$-eigenpairs. Denote
\bqn
T(m,n)&:=&n(m-1)^{n-1},   \\
 E(m,n)&:=&((m-1)^n-1)/(m-2),  \\
G(m,m',n)&:=&((m-1)^n-(m'-1)^n)/(m-m').
\eqn
Recall that Theorem~\ref{theorem2} explains that for tensors  $\cA\in \bC^{[m,n]}$ and $\cB\in \bC^{[m',n]}$, the number of equivalence classes of isolated $\cB$-eigenpairs of $\cA$ is bounded by $T(m,n)$ for $m=m'$ and $G(m,m',n)$ for $m\ne m'$. In particular, Remark~\ref{remark2} states that the number of equivalence classes of isolated eigenpairs of the eigenproblem (\ref{Heig}) is bounded by $T(m,n)$ and
the number of equivalence classes of isolated E-eigenpairs of $\cA$ is bounded by $E(m,n)$. 

\begin{60}
\label{example5-18}
{\rm
In this example, we compare the performance of our \texttt{TenEig}  with \texttt{NSolve} and \texttt{PSOLVE}.  \texttt{teig}, \texttt{eeig}, \texttt{NSolve} and \texttt{PSOLVE} are used to compute eigenpairs (\ref{Heig}) or E-eigenpairs of a generic tensor $\cA \in \bC^{[m,n]}$. We remark the following:

(a) \texttt{teig} is based on Algorithm~\ref{algorithm1}. The polynomial system solved is (\ref{qieig}). 

(b) \texttt{eeig} is based on Algorithm~\ref{algorithm2}. The polynomial system solved is (\ref{polyhyperplane}).

(c)  When computing eigenpairs (\ref{Heig}), \texttt{NSolve} and \texttt{PSOLVE} solve the polynomial system defined by (\ref{qieig}).

(d) When computing E-eigenpairs, \texttt{NSolve} and \texttt{PSOLVE} solve the polynomial system defined by (\ref{evalue}).

The  tensors $\cA$ were generated using $randn(n,\cdots,n)+i*randn(n,\cdots,n)$ in MATLAB. The results are given in Table~\ref{table5-18}. In this table, $N$ denotes the number of equivalence classes of eigenpairs (\ref{Heig}) or E-eigenpairs found by \texttt{teig}, \texttt{eeig}, \texttt{PSOLVE} or \texttt{NSolve}, the reported CPU times are in seconds, ``-'' denotes that no results were returned after 12 hours.

\begin{table}[htbp]
\footnotesize
\begin{center}
\begin{tabular}{|c|c|c|c|c|| c|c|c|c|}
\hline
$(m,n)$ & $T(m,n)$ & Alg & $N$ & time (s) &  $E(m,n)$ & Alg & $N$ & time (s) \\
\hline
\multirow{3}{*}{$(4,5)$} 	&  \multirow{3}{*}{$405$}  & \texttt{teig} & 405 & 15.8 
				&  \multirow{3}{*}{$121$}  & \texttt{eeig} & 121 & 5.4 \\
\cline{3-5} \cline{7-9}
	& & \texttt{PSOLVE} & 404 &   14.0 
	&  & \texttt{PSOLVE} & 121 &   9.5   \\
\cline{3-5} \cline{7-9}
	& & \texttt{NSolve} & 405 &   3136.4 
	&  & \texttt{NSolve} & 121 &   486.6   \\
\hline
\multirow{3}{*}{$(5,5)$}  & \multirow{3}{*}{$1280$}  & \texttt{teig} & 1280 & 73.8 
				& \multirow{3}{*}{$341$}  & \texttt{eeig} &  341 & 22.3 \\
\cline{3-5} \cline{7-9}
	& & \texttt{PSOLVE} & 1280 &   65.5 
	&  & \texttt{PSOLVE} & 341 &   38.6   \\
\cline{3-5} \cline{7-9}
	& & \texttt{NSolve} & - &   - 
	&  & \texttt{NSolve} & 341 &   9264.8   \\
\hline
\multirow{3}{*}{$(5,6)$} 	&  \multirow{3}{*}{$6144$}  & \texttt{teig} & 6144 & 606.5 
				&  \multirow{3}{*}{$1365$}  & \texttt{eeig} &  1365 & 166.5 \\
\cline{3-5} \cline{7-9}
	& & \texttt{PSOLVE} & 6144 &   694.2
	&  & \texttt{PSOLVE} & 1365 &   283.6   \\
\cline{3-5} \cline{7-9}
	& & \texttt{NSolve} & - &   - 
	&  & \texttt{NSolve} & - &   -   \\
\hline
\multirow{3}{*}{$(6,6)$} 	&  \multirow{3}{*}{$18750$}  & \texttt{teig} & 18750 & 3721.3 
				&  \multirow{3}{*}{$3906$}  & \texttt{eeig} & 3906 & 990.2 \\
\cline{3-5} \cline{7-9}
	& & \texttt{PSOLVE} & 18748 &   4636.0
	&  & \texttt{PSOLVE} & 3905 &   1721.0   \\
\cline{3-5} \cline{7-9}
	& & \texttt{NSolve} & - &   - 
	&  & \texttt{NSolve} & - &   -   \\
\hline
\end{tabular}
\end{center}
\caption{Comparison  of \texttt{teig} and \texttt{eeig} with \texttt{PSOLVE} and \texttt{NSolve}}
\label{table5-18}
\end{table}

From table~\ref{table5-18}, we  see that \texttt{teig} and \texttt{eeig} successfully find all equivalence classes of eigenpairs (\ref{Heig}) or E-eigenpairs using reasonable amount of time in all the cases.  \texttt{NSolve} cannot get any results in 12 hours  in some cases (we terminated it after 12 hours).  Although \texttt{PSOLVE} successfully finds all equivalence classes  in many cases, it does miss a few equivalence classes in some cases. We believe that the robustness of \texttt{teig} and \texttt{eeig} is due to their use of the retracing strategies and working in the projective space, as described in Subsection 3.3.  Regarding  the CPU time usage,  \texttt{PSOLVE} is comparable to \texttt{teig} when the number of eigenpairs (\ref{Heig}) $T(m,n)$ is moderate. When $T(m,n)$ gets larger, \texttt{teig} uses less time. This is because the mixed volume computation in \texttt{PSOLVE} takes significantly more time when $T(m,n)$ becomes large.  We also observe that \texttt{eeig} uses less time than \texttt{PSOLVE}. Note that by using the equivalent class structure discussed in Remark \ref{remark4},  \texttt{eeig} traces $E(m,n)$ paths. On the other hand, \texttt{PSOLVE} traces $2E(m,n)$ paths because here it solves the system (\ref{evalue}) from Definition \ref{def1} directly.      
}
\end{60}

\begin{60}
\label{example5-16}
{\rm
 In this example we show the effectiveness and efficiency of \texttt{teig} for finding all equivalence classes of isolated eigenpairs as defined in (\ref{Heig}) of a generic tensor $\cA \in \bC^{[m,n]}$.  Each tensor was generated using $randn(n,\cdots,n)+i*randn(n,\cdots,n)$ in MATLAB. The results are reported in Table~\ref{table5-16}, in which $N$ denotes  the number of equivalence classes of isolated eigenpairs  found by \texttt{teig} and  $T(m,n)$ denotes the bound of the number of equivalence classes of isolated eigenpairs (see Remark \ref{remark2}(a)). 
\begin{table}[htbp]
\begin{center}
\begin{tabular}{|c|c|c|c|| c|c|c|c|}
\hline
$(m,n)$ & $T(m,n)$ & $N$ & time(s) & $(m,n)$ & $T(m,n)$ & $N$ & time(s) \\ \hline
$(3,5)$ & 80 & 80 & 2.4 &  $(3,6)$ & 192 & 192 & 6.8\\ \hline
$(3,7)$ & 448 & 448 & 18.3 & $(3,8)$ & 1024 & 1024& 53.0 \\ \hline
$(3,9)$ & 2304 & 2304 & 145.9 & $(3,10)$ & 5120 & 5120 & 409.2 \\ \hline
$(4,3)$ & 27 & 27 & 0.7 & $(4,4) $ & 108 & 108 & 2.9 \\ \hline
$(4,5)$ & 405 & 405 & 15.8 & $(4,6) $ & 1458 & 1458 & 80.0 \\ \hline
$(4,7)$ & 5103 & 5103 & 385.9 & $(4,8) $ & 17496 & 17496 & 2115.5 \\ \hline
$(5,3)$ & 48 & 48 & 1.2 & $(5,4)$ & 256 & 256 & 8.8 \\ \hline
$(5,5)$ & 1280 & 1280 & 73.8 & $(5,6)$ & 6144 & 6144 & 606.5 \\ \hline
$(5,7)$ & 28672 & 28672 & 5394.2 & $(6,3)$ & 75 & 75 & 2.3 \\ \hline
$(6,4)$ & 500 & 500 & 21.0 & $(6,5)$ & 3125 & 3125 & 287.7 \\ \hline
$(6,6)$ & 18750 & 18750 & 3721.3 & $(7,3)$  & 108 & 108 & 3.6 \\ \hline 
$(7,4)$ & 864 & 864 & 51.5 & $(7,5)$ & 6480 & 6480 & 981.3 \\ \hline
\end{tabular}
\end{center}
\caption{Performance of \texttt{teig} on computing eigenpairs  (\ref{Heig}) of complex random tensors}
\label{table5-16}
\end{table}
}
\end{60}

\begin{60}
\label{example5-17}
{\rm
In this example we show the effectiveness and efficiency of \texttt{eeig} for finding all equivalence classes of isolated E-eigenpairs of a generic tensor $\cA\in \bC^{[m,n]}$.  Each generic tensor was generated using $randn(n,\cdots,n)+i*randn(n,\cdots,n)$ in MATLAB. The results are reported in Table~\ref{table5-17}, in which $N$ denotes  the number of equivalence classes of E-eigenpairs found by \texttt{eeig} and  $E(m,n)$ denotes the bound of the number of equivalence classes of isolated E-eigenpairs (see Remark \ref{remark2}(b)). 
\begin{table}[htbp]
\begin{center}
\begin{tabular}{|c|c|c|c|| c|c|c|c|}
\hline
$(m,n)$ & $E(m,n)$ & $N$ & time(s) & $(m,n)$ & $E(m,n)$ & $N$ & time(s) \\ \hline
$(3,5)$ & 31 & 31 & 1.4 & $(3,6)$ & 63 & 63 & 3.1 \\ \hline
$(3,7)$ & 127 & 127 & 7.5 & $(3,8)$ & 255 & 255 & 20.3 \\ \hline
$(3,9)$ & 511 & 511 & 48.5 & $(3,10)$ & 1023 & 1023 & 133.9 \\ \hline
$(4,3)$ & 13 & 13 & 0.4 & $(4,4)$ & 40 & 40 & 1.7 \\ \hline
$(4,5)$ & 121 & 121 & 5.4 & $(4,6) $ & 364 & 364 & 26.9 \\ \hline
$(4,7)$ & 1093 & 1093 & 119.5 & $(4,8) $ & 3280 & 3280 & 555.8 \\ \hline
$(5,3)$ & 21 & 21 & 0.7 & $(5,4)$ & 85 & 85 & 4.2 \\ \hline
$(5,5)$ & 341 & 341 & 22.3 & $(5,6)$ & 1365 & 1365 & 166.5 \\ \hline
$(5,7)$ & 5461 & 5461 & 1330.7 & $(6,3)$ & 31 & 31 & 1.2  \\ \hline
$(6,4)$ & 156 & 156 & 9.5 & $(6,5)$ & 781 & 781 & 100.4 \\ \hline
$(6,6)$ & 3906 & 3906 & 990.2 & $(7,3)$ & 43 & 43 & 1.9 \\ \hline
$(7,4)$ & 259 & 259 & 21.3 & $(7,5)$ & 1555 & 1555 & 245.0 \\ \hline
\end{tabular}
\end{center}
\caption{Performance of \texttt{eeig} on computing E-eigenpairs of complex random tensors}
\label{table5-17}
\end{table}
}
\end{60}

According to \cite{Qi05}, \cite{CQZ13}, and \cite{CS13}, for a randomly generated  tensor $\cA \in \bC^{[m,n]}$, it has $T(m,n)$  equivalence classes of eigenpairs  (\ref{Heig}) and $E(m,n)$ equivalence classes of E-eigenpairs. Moreover, its eigenpairs and E-eigenpairs are isolated. 
From Tables \ref{table5-16} and \ref{table5-17}, we observe that  \texttt{teig} and \texttt{eeig} can find all equivalence classes of eigenpairs  (\ref{Heig}) and E-eigenpairs of such an tensor in the examples we tested.

\begin{60}
\label{example5-20}
{\rm
In this example we show the effectiveness and efficiency of \texttt{teig} and \texttt{teneig} for finding all equivalence classes of isolated $\cB$-eigenpairs of a tensor $\cA$, where $\cA\in \bC^{[m,n]}, \cB \in \bC^{[m',n]}$ are generic tensors. Each generic tensor was generated using $randn(n,\cdots,n)+i*randn(n,\cdots,n)$ in MATLAB. The results are reported in Table~\ref{table5-20}, in which $N$ denotes  the number of equivalence classes of eigenpairs found by \texttt{teig} or \texttt{teneig}, $T(m,n)$ denotes the bound of the number of equivalence classes of isolated $\cB$-eigenpairs for $m=m'$, and $G(m,m',n)$ denotes the bound of the number of equivalence classes of isolated $\cB$-eigenpairs for $m\ne m'$ (see Theorem~\ref{theorem2}). 
\begin{table}[htbp]
\begin{center}
\begin{tabular}{|c|c|c|c|| c|c|c|c|}
\hline
\multicolumn{4}{|c||}{teig ($m = m'$)} & \multicolumn{4}{|c|}{teneig ($m \ne m'$)} \\ \hline
$(m,n)$ & $T(m,n)$ & $N$ & time(s) & $(m,m',n)$ & $G(m,m',n)$ & $N$ & time(s) \\ \hline
$(3,7)$ & 448 & 448 & 23.7 	& $(3,2,7)$ & 127 & 127 & 10.3 \\ \hline
$(3,8)$ & 1024 & 1024 & 68.3 	& $(3,4,6)$ & 665 & 665 & 68.1 \\ \hline
$(3,9)$ & 2304 & 2304 & 210.3 	& $(3,5,5)$ & 496 & 496 & 49.5 \\ \hline
$(4,5)$ & 405 & 405 & 20.8 	& $(4,2,6)$ & 364 & 364 & 28.9 \\ \hline
$(4,6)$ & 1458 & 1458 & 110.4 	& $(4,3,5)$ & 211 & 211 & 13.2 \\ \hline
$(4,7)$ & 5103 & 5103& 737.5 	& $(4,5,4)$ & 175 & 175 & 9.5 \\ \hline
$(5,5)$ & 1280 & 1280 & 97.9 	& $(5,4,5)$ & 781 & 781 & 83.9 \\ \hline
$(5,6)$ & 6144 & 6144 & 623.6 	& $(5,6,3)$ & 61 & 61 & 2.6 \\ \hline
$(6,4)$ & 500 & 500 & 29.9 	& $(6,5,4)$ & 369 & 369 & 30.7  \\ \hline
$(6,5)$ & 3125 & 3125 & 449.4 	& $(6,7,3)$ & 91 & 91 & 6.0 \\ \hline
$(7,3)$ & 108 & 108 & 4.4 	& $(7,6,4)$ & 671 & 671 & 77.1 \\ \hline
$(7,4)$ & 864 & 864 & 77.6 	& $(7,8,3)$ & 127 & 127 & 9.4 \\ \hline
\end{tabular}
\end{center}
\caption{Performance of \texttt{teig} and \texttt{teneig} on computing generalized eigenpairs of complex random tensors}
\label{table5-20}
\end{table}
}
\end{60}

From Table~\ref{table5-20} we see that our \texttt{teig} and \texttt{teneig} find all equivalence classes of isolated $\cB$-eigenpairs of $\cA$ for the generic tensors $\cA$ and $\cB$ we tested in a reasonable amount of time.

\subsection{Examples for Computing Real Eigenpairs}

In this subsection, numerical examples are provided to illustrate the effectiveness and efficiency of \texttt{zeig} or \texttt{heig} for computing real Z-eigenpairs or H-eigenpairs of a tensor $\cA \in \bR^{[m,n]}$.  
By Definition~\ref{def1},  $(\lambda,x)$ is a Z-eigenpair if and only if $((-1)^{m-2}\lambda,-x)$ is a Z-eigenpair, and $(\lambda,x)$ is an H-eigenpair if and only if $(\lambda,tx)$ is an H-eigenpair for any nonzero
$t \in \bR$. Only one representative from each equivalence class of eigenpairs will be listed in our examples.  The notation $\lambda^{(l)}$ is used to denote $l$ eigenvectors counting multiplicities are found for the eigenvalue $\lambda$. In the following tables, the multiplicity of an eigenpair means the multiplicity of this eigenpair as a zero of the corresponding defining polynomial system. For the sake of conciseness, the polynomial system resulted from the tensor eigenvalue problem will be omitted.  

\begin{60}
\label{example5-1}
{\rm
Consider the symmetric tensor $\cA \in \bR^{[6,3]}$  whose corresponding polynomial form is the Motzkin polynomial:  
$$
\cA x^6 =x_{3}^6+x_{1}^4x_{2}^2+x_{1}^2x_{2}^4-3x_{1}^2x_{2}^2x_{3}^2.
$$
In Example 5.9 of \cite{CS13}, it states that this tensor has 25 equivalence classes of Z-eigenpairs. Using \texttt{zeig} exactly 25 equivalence classes of Z-eigenpairs are found as shown in Table~\ref{table5-1}, which confirms the results of \cite{CS13}. \texttt{zeig} takes about 0.9 seconds
to carry out the entire computation.

\begin{table}[htbp]
\begin{center}
\begin{tabular}{ |c|c|c|c|c|c|c|c| }
\hline
$\lambda$ & \multicolumn{3}{|c|}{$0^{(14)}$} & \multicolumn{2}{|c|}{$0.0156^{(8)}$}  & $0.2500^{(2)}$ & $1$\\ \hline
$x_1$ & $ 0.5774$ & $1$ & $0$ & $ 0.8253$ & $ 0.2623$ & $ 0.7071 $ & $0$\\ \hline
$x_2$ & $\pm 0.5774$ & $0$ & $1$ & $\pm 0.2623$ & $\pm 0.8253$ & $\pm 0.7071 $ & $0$\\ \hline
$x_3$ & $\pm 0.5774$ & $0$ & $0$ & $\pm 0.5000$ & $\pm 0.5000$ & $0 $ & $1$\\ \hline
multiplicity & 1 & 5 & 5 & 1 & 1 & 1 & 1\\ \hline
\end{tabular}
\caption{Z-eigenpairs of the tensor in Example~\ref{example5-1}}
\label{table5-1}
\end{center}
\end{table}

All the H-eigenpairs found by Example 4.10 of \cite{CDN14} are also found by \texttt{heig} as shown in Table~\ref{table5-21}. \texttt{heig} takes about 1.7 seconds to carry out the entire computation.

\begin{table}[htbp]
\begin{center}
\begin{tabular}{ |c|c|c|c|c|c|c|c|}
\hline
$\lambda$ & \multicolumn{3}{|c|}{$0^{(14)}$} & \multicolumn{2}{|c|}{$0.0555^{(8)}$}  & \multicolumn{2}{|c|}{$1^{(15)}$} \\ \hline
$x_1$ & $ \pm 1$ & $1$ & $0$ & $ \pm 0.4568$ & $ 1$ & $ 0 $ & $1$\\ \hline
$x_2$ & $1$ & $0$ & $1$ & $1$ & $\pm 0.4568$ & $0 $ & $\pm 1$ \\ \hline
$x_3$ & $\pm 1$ & $0$ & $0$ & $\pm 0.6856$ & $\pm 0.6856$ & $1 $ & $1$\\ \hline
multiplicity & 1 & 5 & 5 & 1 & 1 & 13 & 1 \\ \hline
\end{tabular}
\caption{H-eigenpairs of the tensor in Example~\ref{example5-1}}
\label{table5-21}
\end{center}
\end{table}
}
\end{60}

So far the only available method that can find all real eigenvalues of a symmetric tensor is Algorithm 3.6 in \cite{CDN14}.  In the next two examples,  we report our experiments on \texttt{zeig} for computing all Z-eigenvalues, using examples taken from \cite{CDN14}.

\begin{60}
\label{examples}
{\rm
We use our \texttt{zeig} to compute the Z-eigenpairs of 12 symmetric tensors from \cite{CDN14}. The test problems and numerical results are given in the Appendix.  From the numerical results we see that \texttt{zeig} finds all the Z-eigenvalues found by Algorithm 3.6 of \cite{CDN14} on this set of test problems. We now summarize the CPU time (in seconds) used by \texttt{zeig}  in Table~\ref{table-examples}.  Since the computer used in \cite{CDN14} is different from the computer used in this paper,  the CPU time used by Algorithm 3.6 in \cite{CDN14} is not reported here, but we refer to \cite{CDN14}.

\begin{table}[htbp]
	\begin{center}
		\begin{tabular}{|c|c|c|c|c|c|c|c|c|c|c|c|c|}
			\hline
			Problem & 1 & 2 & 3 & 4 & 5 & 6 & 7 & 8 & 9 & 10 & 11 & 12 \\
			\hline
			time(s) & 0.3 & 4.0 & 0.3--0.4 & 0.1 & 0.6 & 1.8 & 15.7 & 6.1 & 0.3 & 6.3 & 27.3 & 4.5 \\
			\hline
		\end{tabular}
	\end{center}
	\caption{CPU time used by \texttt{zeig} for computing the Z-eigenvalues of 12 symmetric tensors from \cite{CDN14}}
	\label{table-examples}
\end{table}

}
\end{60}

\begin{60}
\label{example5-19}
{\rm
Consider the symmetric tensor $\cA \in \bR^{[4,n]}$ (Example 4.16 in \cite{CDN14}) with the polynomial form
\bqn
\cA x^4 &=& (x_1-x_2)^4 + \cdots + (x_1-x_n)^4 + (x_2-x_3)^4 + \cdots + (x_2-x_n)^4 \\
&& + \cdots + (x_{n-1}-x_n)^4.
\eqn
For different $n$, all the  Z-eigenvalues found by Algorithm 3.6 in \cite{CDN14} are also found by \texttt{zeig}, which are given in Table \ref{table5-19}. 
We remark  that when $n=8,9,10$,  \texttt{zeig} can find all the Z-eigenvalues in a reasonable amount of time, but \cite{CDN14} reports that Algorithm 3.6 can only find the first three largest Z-eigenvalues. The CPU time used by \texttt{zeig} is reported in Table~\ref{table5-19}. Since different computers were used, we refer to \cite{CDN14} for the CPU time used by Algorithm 3.6 of \cite{CDN14}.  For the sake of conciseness, the corresponding Z-eigenvectors are not displayed.

\begin{table}[htbp]
	\begin{center}
		\begin{tabular}{ |c|ccccc|c|}
			\hline
			$n$ & \multicolumn{5}{ c| }{$\lambda$} & time(s) \\ \hline
			4 & 0.0000 & 4.0000 & 5.0000 & 5.3333 & & 1.7   \\ \hline 
			5 & 0.0000, & 4.1667, & 4.2500, & 5.5000, & 6.2500  & 5.4  \\ \hline 
			6 & 0.0000, & 4.0000, &4.5000, & 6.0000, & 7.2000 & 15.5   \\ \hline 
			\multirow{2}{*}{7} & 0.0000, & 4.0833, & 4.1667, & 4.7500, & 4.8846, & \multirow{2}{*}{58.3}  \\
			& 4.9000, & 6.5000, & 8.1667 & & &  \\ \hline 
			\multirow{2}{*}{$8^*$}  & 0.0000,& 4.0000, & 4.2667,& 4.2727, & 4.3333,& \multirow{2}{*}{244.1}   \\ 
			& 5.0000,& 5.2609, &5.3333, & 7.0000, & 9.1429 &   \\ \hline 
			\multirow{2}{*}{$9^*$}  & 0.0000,& 4.0500, & 4.1250,& 4.5000, & 5.2500,& \multirow{2}{*}{788.0} \\ 
			& 5.6250,& 5.7857, & 7.5000, & 10.1250 &   &   \\ \hline 
			\multirow{3}{*}{$10^*$}  & 0.0000,& 4.0000, & 4.1667,& 4.1818, & 4.2500,& \multirow{3}{*}{2665.6} 
			\\  
			& 4.6667,& 4.7500, & 4.7593, & 4.7619, & 5.5000,  &  \\ 
			& 5.9808, & 6.2500, & 8.0000, & 11.1111 &      &  \\ \hline 
		\end{tabular}
		\caption{CPU time of \texttt{zeig} for computing Z-eigenvalues of the tensor in Example~\ref{example5-19} ($*$ denotes that Algorithm 3.6 (\cite{CDN14}) only finds the first three largest Z-eigenvalues.) }
		\label{table5-19}
	\end{center}
\end{table}

}
\end{60}

\ \\
{\bf Acknowledgments.} We would like to thank Professor T. Y. Li,  the two anonymous referees, and the associate editor for their insightful and constructive comments and suggestions, which significantly improved both the content and presentation of the paper.


\newpage

\begin{center}
{\Large \bf Appendix}  
\end{center}

\begin{80}
\label{example5-3}
{\rm
Consider the symmetric tensor $\cA \in \bR^{[4,3]}$ (Example 4.1 in \cite{CDN14}, see also \cite{Qi05}) with the polynomial form
$$
\cA x^4 = x_{1}^4 + 2x_{2}^4 + 3 x_{3}^4.
$$
Using \texttt{zeig} all the Z-eigenpairs found in \cite{CDN14} are obtained (see Table~\ref{table5-3}). \texttt{zeig} takes about 0.3 seconds to carry out the entire computation. 
\begin{table}[htbp]
\begin{center}
\begin{tabular}{ |c|c|c|c|c|c|c|c| }
\hline
$\lambda$ & $0.5455^{(4)}$ & $0.6667^{(2)}$  & $0.7500^{(2)}$ & $1$ & $1.2^{(2)}$ & $2$ & $3$\\ \hline
$x_1$ & $ 0.7385$ & $ 0.8165$ & $ 0.8660$ & $ 1$ & $0$ & $0$ & $0$ \\ \hline
$x_2$ & $\pm 0.5222$ & $\pm 0.5774$ & $0$ & $0$ & $ 0.7746$ & $1$ & $0$ \\ \hline
$x_3$ & $\pm 0.4264$ & $0$ & $\pm 0.5000$  & $0$ & $\pm 0.6325$ & $0$ & $1$\\ \hline
multiplicity & 1 & 1 & 1  & 1 & 1 & 1 & 1\\ \hline
\end{tabular}
\caption{Z-eigenpairs of the tensor in Problem~\ref{example5-3}}
\label{table5-3}
\end{center}
\end{table}
}
\end{80}

\begin{80}
\label{example5-4}
{\rm
For the diagonal tensor $\cD \in \bR^{[5,4]}$ (Example 4.2 in \cite{CDN14}) such that $\cD x^5 = x_1^5 + 2x_2^5 - 3x_3^5 - 4x_4^5$. Consider the symmetric tensor $\cA \in \bR^{[5,4]}$ such that $\cA x^5 = \cD(Qx)^5$ where
$$
Q = (I-2w_1w_1^T)(I-2w_2w_2^T)(I-2w_3w_3^T)
$$
and $w_1,w_2,w_3$ are randomly generated unit vectors. All the 30 Z-eigenpairs found in \cite{CDN14} are also found by using \texttt{zeig}. The 15 nonnegative Z-eigenvalues are listed below
\[
\begin{array}{cccccccc} 
0.2518, &  0.3261,  &  0.3466, &  0.3887, &   0.4805, &   0.5402,  &  0.5550, &   0.6057, \\
0.8543, &  0.9611,  & 1.0000,  &  1.2163, &   2.0000,  &  3.0000, &   4.0000. &
\end{array}
\]
For conciseness,  the corresponding Z-eigenvectors are not displayed here. \texttt{zeig} takes about 4.0 seconds to do the entire computation.
}
\end{80}

\begin{80}
\label{example5-5}
{\rm
Consider the symmetric tensor $\cA \in \bR^{[4,3]}$ (Example 4.3 in \cite{CDN14}, see also \cite{Qi05}) with the polynomial form
$$
\cA x^4 = 2x_1^4 + 3x_2^4 + 5x_3^4 + 4ax_1^2x_2x_3,
$$
where $a$ is a parameter. All Z-eigenvalues  found in \cite{CDN14} are also found by \texttt{zeig} for different values of $a$, as shown in Table~\ref{table5-5}.
\begin{table}[htbp]
\begin{center}
\begin{tabular}{ |c|c|c|}
\hline
$a$ & $\lambda$ & time(s) \\ \hline
0 & $0.9677^{(4)}$, $1.2000^{(2)}$, $1.4286^{(2)}$, $1.8750^{(2)}$, $2$, $3$, $5$ & 0.4\\ \hline
0.25 & $0.8464^{(2)}$, $1.0881^{(2)}$, $1.2150^{(2)}$, $1.4412^{(2)}$, $1.8750^{(2)}$, $2$, $3$, $5$ & 0.4 \\ \hline
0.5 & $0.7243^{(2)}$, $1.2069^{(2)}$, $1.2593^{(2)}$, $1.4783^{(2)}$, $1.8750^{(2)}$, $2$, $3$, $5$ & 0.4\\ \hline
1 & $0.4787^{(2)}$, $1.6133^{(2)}$, $1.8750^{(2)}$, $2$, $3$, $5$ & 0.3 \\ \hline
3 & $-0.5126^{(2)}$, $1.8750^{(2)}$, $2$, $2.2147^{(2)}$, $3$, $5$ & 0.3\\ \hline
\end{tabular}
\caption{Z-eigenvalues of the tensor in Problem~\ref{example5-5}}
\label{table5-5}
\end{center}
\end{table}
For conciseness, the corresponding Z-eigenvectors are not displayed here. The CPU time used by \texttt{zeig} for each $a$ is also reported in the table. 
}
\end{80}

\begin{80}
\label{example5-7}
{\rm
Consider the symmetric tensor $\cA \in \bR^{[4,2]}$ (Example 4.4 in \cite{CDN14}, see also \cite{Qi05}) with the polynomial form
$$
\cA x^4 = 3x_1^4 + x_2^4 + 6ax_1^2x_2^2,
$$
where $a$ is a parameter. All Z-eigenvalues found in \cite{CDN14} are also found by \texttt{zeig} for different values of $a$, which are  listed in Table~\ref{table5-7}.  The CPU time used 
by \texttt{zeig} for each $a$ is also given in the table. 
\begin{table}[htbp]
\begin{center}
\begin{tabular}{ |c|c|c|}
\hline
$a$ & $\lambda$ & time(s) \\ \hline
$-1$ & $-0.6000^{(2)}$, $1$, $3$ & 0.1 \\ \hline
$0$ &  $0.7500^{(2)}$, $1$, $3$ & 0.1 \\ \hline
$0.25$ & $0.9750^{(2)}$, $1$, $3$ & 0.1 \\ \hline
$0.5$ & $1$, $3$ & 0.1 \\ \hline
$2$ & $1$, $3$, $4.1250^{(2)}$ & 0.1 \\ \hline
$3$ & $1$, $3$, $5.5714^{(2)}$ & 0.1 \\ \hline
\end{tabular}
\caption{Z-eigenvalues of the tensor in Problem~\ref{example5-7}}
\label{table5-7}
\end{center}
\end{table}
For conciseness, the corresponding Z-eigenvectors are not displayed here.
}
\end{80}

\begin{80}
\label{example5-8}
{\rm
Consider the symmetric tensor $\cA \in \bR^{[4,3]}$ (Example 4.5 in \cite{CDN14}, see also \cite{KM11} or \cite{NW13}) such that
\[
\begin{array}{c} 
A_{1111}= 0.2883, A_{1112}=  -0.0031, A_{1113}=0.1973, A_{1122}= -0.2485,  \\
A_{1123}= -0.2939, A_{1133}=0.3847, A_{1222}=0.2972, A_{1223}=0.1862, \\
A_{1233}=0.0919, A_{1333}=-0.3619, A_{2222}=0.1241, A_{2223}=-0.3420,\\
A_{2233}=0.2127, A_{2333}=0.2727, A_{3333}=-0.3054.
\end{array}
\]
All the Z-eigenpairs found in \cite {CDN14} are also found by \texttt{zeig}, as given in Table \ref{table5-8}. 
\texttt{zeig} takes about 0.6 seconds to do the entire computation.
\begin{table}[htbp]
\scriptsize
\begin{center}
\begin{tabular}{ |c|c|c|c|c|c|c|c|c|c|c|c|}
\hline
$\lambda$ & -1.0954 & -0.5629 & -0.0451  &  0.1735  &  0.2433 & 0.2628 & 0.2682 & 0.3633 & 0.5105  &  0.8169 & 0.8893\\ \hline
$x_1$ & -0.5915  & -0.1762 & 0.7797 & 0.3357  & -0.9895   &   -0.1318    &   0.6099   &   0.2676  &  -0.3598 & -0.8412    &   -0.6672 \\ \hline
$x_2$ & 0.7467   &   0.1796   &   0.6135    &   0.9073   & -0.0947    &    0.4425   &   0.4362   &   0.6447  &  0.7780   &    0.2635  &   -0.2471 \\ \hline
$x_3$ & 0.3043  &  -0.9678    &   0.1250   &   0.2531    &   0.1088  &   0.8870  &   0.6616  &   0.7160  & -0.5150    &   -0.4722  &   0.7027 \\ \hline
\end{tabular}
\caption{Z-eigenpairs of the tensor in Problem~\ref{example5-8}}
\label{table5-8}
\end{center}
\end{table}
}
\end{80}

\begin{80}
\label{example5-9}
{\rm
Consider the symmetric tensor $\cA \in \bR^{[3,6]}$ (Example 4.6 in \cite{CDN14}, see also \cite{QWW09}) such that
$A_{iii}= i$ for $i = 1,\dots,6$ and $A_{i,i,i+1}=10$ for $i=1,\dots,5$ and zero otherwise. All the Z-eigenvalues found in \cite{CDN14} are also found by \texttt{zeig}.  The 19 nonnegative Z-eigenvalues are listed below:  
\[
\begin{array}{cccccccc} 
3.9992  &  4.0225  &  4.2464  &  4.3358  &  5.1402  & 5.4817  &  5.5218  &  5.5668  \\
 5.5674  &   6.0000  &  7.2165  &  8.1889  &  8.5979 & 8.6596  &  8.7347 & 10.9711  \\
 15.4298  & 15.4552  & 16.2345 &  & & &
\end{array}
\]
For conciseness, the corresponding Z-eigenvectors are not displayed here. \texttt{zeig} takes about 1.8 seconds to carry out the entire computation.
}
\end{80}

\begin{80}
\label{example5-10}
{\rm
Consider the symmetric tensor $\cA \in \bR^{[4,6]}$ (Example 4.7 in \cite{CDN14}, see also \cite{LQY11}) with the polynomial form
\bqn
-\cA x^4 &=& (x_1-x_2)^4 + (x_1-x_3)^4 + (x_1-x_4)^4 + (x_1-x_5)^4 + (x_1-x_6)^4\\
&& + (x_2-x_3)^4  + (x_2-x_4)^4  + (x_2-x_5)^4  + (x_2-x_6)^4 \\
&&  + (x_3-x_4)^4 + (x_3-x_5)^4 + (x_3-x_6)^4 \\
&&  + (x_4-x_5)^4 + (x_4-x_6)^4 + (x_5-x_6)^4.
\eqn
All the 5 Z-eigenvalues found in \cite{CDN14} are also found by \texttt{zeig}, which are given in Table~\ref{table5-10}.  
\begin{table}[htbp]
\begin{center}
\begin{tabular}{ |c|c|c|}
\hline
$\lambda$ & $x^T$ & multiplicity \\ \hline
$-7.2000^{(6)}$ & $(0.1826,0.1826,0.1826,0.1826,0.1826,-0.9129)$ & 1\\ \hline
$-6.0000^{(15)}$ & $(0.7071,0,0,0,0,-0.7071)$ & 1\\ \hline
$-4.5000^{(\star)}$ & $(0.5774, 0.5774, -0.2887, -0.2887, -0.2887, -0.2887)$  & -\\ \hline
$-4.0000^{(10)}$ & $(0.4082,0.4082,0.4082,-0.4082,-0.4082,-0.4082)$ & 1\\ \hline
$0^{(\star)}$ & $(0.4082, 0.4082, 0.4082, 0.4082, 0.4082, 0.4082)$  & - \\ \hline
\end{tabular}
\caption{Z-eigenpairs of the tensor in Problem~\ref{example5-10}}
\label{table5-10}
\end{center}
\end{table}
As pointed out in \cite{CDN14}, every permutation of a Z-eigenvector is also a Z-eigenvector. Only the Z-eigenvector with $x_1\ge x_2 \ge \cdots \ge x_6$ corresponding to one Z-eigenvalue is listed. We remark that the Z-eigenpairs corresponding to Z-eigenvalues $0$ and $-4.5$ are in a positive dimensional solution component of the corresponding polynomial system. Therefore, there are  infinitely many Z-eigenvectors associated with 0 and $-4.5$. \texttt{zeig} finds 484 Z-eigenvectors associated with $0$ and 180 Z-eigenvectors associated with $-4.5$. Only one of these Z-eigenvectors for each case is listed in Table~\ref{table5-11}. \texttt{zeig} takes about 15.7 seconds to do the entire computation.
}
\end{80}

\begin{80}
\label{example5-11}
{\rm
Consider the symmetric tensor $\cA \in \bR^{[4,5]}$ (Example 4.8 in \cite{CDN14}, see also \cite{XC13}) with the polynomial form
$$
\cA x^4 = (x_1+x_2+x_3+x_4)^4 +  (x_2+x_3+x_4+x_5)^4 .
$$
All the 3 Z-eigenvalues found in \cite{CDN14} are also found by \texttt{zeig}, which are shown in Table~\ref{table5-11}. 
\begin{table}[htbp]
\begin{center}
\begin{tabular}{ |c|c|c|}
\hline
$\lambda$ & $x^T$ & multiplicity \\ \hline
$0^{(\star)}$ & $( 0.3870, -0.1537, 0.4532, -0.6866, 0.3870)$  & - \\ \hline
0.5000 & $(0.7071,0,0,0,-0.7071)$ & 1 \\ \hline
24.5000 & $(0.2673,0.5345,0.5345,0.5345,0.2673)$ & 1 \\ \hline
\end{tabular}
\caption{Z-eigenpairs of the tensor in Problem~\ref{example5-11}}
\label{table5-11}
\end{center}
\end{table}

We remark that the Z-eigenpairs corresponding to Z-eigenvalue 0 are in a positive dimensional solution component of the corresponding polynomial system. Thus, there are infinitely many Z-eigenvectors associated with Z-eigenvalue 0. \texttt{zeig} finds 234  of them. Only one of them is listed in Table~\ref{table5-11}.  \texttt{zeig} uses about 6.1 seconds to do the entire computation.
}
\end{80}

\begin{80}
\label{example5-12}
{\rm
Consider the symmetric tensor $\cA \in \bR^{[3,3]}$ (Example 4.9 in \cite{CDN14}, see also \cite{CS13}) with the polynomial form
$$
\cA x^3 = 2x_1^3 + 3x_1x_2^2 + 3x_1x_3^2.
$$
We remark that the Z-eigenpairs corresponding to Z-eigenvalue 2 are in a positive dimensional solution component of the corresponding polynomial system. Thus, there are infinitely many Z-eigenvectors associated with Z-eigenvalue 0.
\texttt{zeig} finds 7  of them. Only one of them is listed in Table~\ref{table5-12}. \texttt{zeig} uses about 0.3 seconds to do the entire computation.
\begin{table}[htbp]
\begin{center}
\begin{tabular}{ |c|c|c|}
\hline
$\lambda$ & $x^T$ & multiplicity \\ \hline
$2^{(\star)}$ & $(1,0,0)$ & - \\ \hline
\end{tabular}
\caption{Z-eigenpairs of the tensor in Problem~\ref{example5-12}}
\label{table5-12}
\end{center}
\end{table}
}
\end{80}

\begin{80}
\label{example5-13}
{\rm
Consider the tensor $\cA \in \bR^{[4,n]}$ (Example 4.12 in \cite{CDN14}, see also \cite{NW13}) such that
$$
\cA_{i_1,\dots,i_4} = \sin (i_1+i_2+i_3+i_4).
$$
When $n=5$, all the 5 Z-eigenvalues found in \cite{CDN14} are also found by \texttt{zeig}, which are given in Table~\ref{table5-13}. 
\begin{table}[htbp]
\begin{center}
\begin{tabular}{ |c|c|c|}
\hline
$\lambda$ & $x^T$ & multiplicity \\ \hline
-8.8463 & $(0.5809, 0.3563,-0.1959,-0.5680,-0.4179)$ & 1 \\ \hline
-3.9204 & $(-0.1785,0.4847,0.7023,0.2742,-0.4060)$ & 1 \\ \hline
$0^{(\star)}$ & $(-0.5213, 0.3748, -0.6608, 0.1824, -0.3433)$ & - \\ \hline
4.6408 & $(0.5055,-0.1228,-0.6382,-0.5669,0.0256)$ & 1 \\ \hline
7.2595 & $(0.2686,0.6150,0.3959,-0.1872, -0.5982)$ & 1 \\ \hline
\end{tabular}
\caption{Z-eigenpairs of the tensor in Problem~\ref{example5-13}}
\label{table5-13}
\end{center}
\end{table}

We remark that the Z-eigenpairs corresponding to Z-eigenvalue $0$ are in a positive dimensional solution component of the corresponding polynomial system. Thus, there are  infinitely many Z-eigenvectors associated with $0$. \texttt{zeig} finds 234  of them. Only one of them is listed in Table~\ref{table5-13}. \texttt{zeig} takes about 6.3 seconds to carry out the entire computation. 
}
\end{80}

\begin{80}
\label{example5-14}
{\rm
Consider the tensor $\cA \in \bR^{[4,n]}$ (Example 4.13 in \cite{CDN14}) such that
$$
\cA_{i_1,\dots,i_4} = \tan (i_1)  + \cdots + \tan (i_4) .
$$
When $n=6$, all the 3 Z-eigenvalues found in \cite{CDN14} are also found by \texttt{zeig}, which are given in Table~\ref{table5-14}. 
\begin{table}[htbp]
\begin{center}
\begin{tabular}{ |c|c|c|}
\hline
$\lambda$ & $x^T$ & multiplicity \\ \hline
-133.2871 & $(0.1936,0.5222,0.3429,0.2287,0.6272,0.3559)$ & 1 \\ \hline
$0^{(\star)}$ & $(-0.5840, -0.3454, 0.1784, 0.6773, 0.1892, -0.1156)$ & - \\ \hline
45.5045 & $(0.6281,0.0717,0.3754,0.5687,-0.1060,0.3533)$ & 1 \\ \hline
\end{tabular}
\caption{Z-eigenpairs of the tensor in Problem~\ref{example5-14}}
\label{table5-14}
\end{center}
\end{table}

We remark that the Z-eigenpairs corresponding to Z-eigenvalue $0$ are in a positive dimensional solution component of the corresponding polynomial system. Thus, there are infinitely many Z-eigenvectors associated with 0. \texttt{zeig} finds 724 of them. Only one of them is listed in Table~\ref{table5-14}. It takes \texttt{zeig} about 27.3 seconds to carry out the entire computation.
}
\end{80}

\begin{80}
\label{example5-15}
{\rm
Consider the tensor $\cA \in \bR^{[5,n]}$ (Example 4.14 in \cite{CDN14}) such that
$$
\cA_{i_1,\dots,i_5} = \ln (i_1)  + \cdots + \ln (i_5) .
$$
For $n=4$, all the 3 Z-eigenvalues found in \cite{CDN14} are also found by \texttt{zeig}, which are shown in Table~\ref{table5-15}. 
\begin{table}[htbp]
\begin{center}
\begin{tabular}{ |c|c|c|}
\hline
$\lambda$ & $x^T$ & multiplicity \\ \hline
$0^{(\star)}$ & $(-0.4304, 0.8139, 0.0069, -0.3903)$ & - \\ \hline
0.7074 & $(-0.9054,-0.3082,0.0411,0.2890)$ & 1\\ \hline
132.3070 & $(0.4040,0.4844,0.5319,0.5657)$ & 1 \\ \hline
\end{tabular}
\caption{Z-eigenpairs of the tensor in Problem~\ref{example5-15}}
\label{table5-15}
\end{center}
\end{table}

We remark that the Z-eigenpairs corresponding to Z-eigenvalue 0 are in a positive dimensional solution component of the corresponding polynomial system. Thus, there are infinitely many Z-eigenvectors associated with 0. \texttt{zeig} finds 166 of them. Only one of them is listed in Table~\ref{table5-14}. The entire computation takes \texttt{zeig} about 4.5 seconds.
}
\end{80}

\end{document}